\documentclass{amsart}
\usepackage{amssymb}
\usepackage{graphics}
\usepackage{latexsym}
\usepackage{amsmath}
\usepackage{amssymb,amsthm,amsfonts}
\usepackage{amscd}
\usepackage[arrow, matrix, curve]{xy}
\usepackage{syntonly}
\title[Higgs bundles over Shimura Curve in char $p$]{Higgs bundles over the good reduction of a quaternionic Shimura
curve}

\author[M. Sheng]{Mao Sheng}

\address{Institut f\"{u}r  Mathematik, Johannes Gutenberg-Universit\"{a}t
Mainz, Mainz, 55099, Germany, Fax: +49 6131 39 21295, Tel: +49 6131
39 27345}

\email{sheng@uni-mainz.de}
\author[J.J. Zhang]{Jiajin Zhang}

\address{Department of Mathematics, Sichuan University, Chengdu, 610000, P.R. China}

\email{jjzhang@scu.edu.cn}

\author[K. Zuo]{Kang Zuo}

\address{Institut f\"{u}r Mathematik, Johannes Gutenberg-Universit\"{a}t
Mainz, Mainz, 55099, Germany }

\email{zuok@uni-mainz.de}

\begin{document}
%%%%%%%%%%%%%%%%%%%% Text italic %%%%%%%%%%%%%%%%%%%%%%%%%%%%
\theoremstyle{plain}
\newcommand{\bbA}{\mathbb{A}}
\newcommand{\bbB}{\mathbb{B}}
\newcommand{\bbC}{\mathbb{C}}
\newcommand{\bbD}{\mathbb{D}}
\newcommand{\bbE}{\mathbb{E}}
\newcommand{\bbF}{\mathbb{F}}
\newcommand{\bbG}{\mathbb{G}}
\newcommand{\bbH}{\mathbb{H}}

\newcommand{\bbK}{\mathbb{K}}
\newcommand{\bbL}{\mathbb{L}}
\newcommand{\bbM}{\mathbb{M}}
\newcommand{\bbN}{\mathbb{N}}
\newcommand{\bbO}{\mathbb{O}}
\newcommand{\bbP}{\mathbb{P}}
\newcommand{\bbQ}{\mathbb{Q}}
\newcommand{\bbR}{\mathbb{R}}
\newcommand{\bbS}{\mathbb{S}}
\newcommand{\bbT}{\mathbb{T}}
\newcommand{\bbU}{\mathbb{U}}
\newcommand{\bbV}{\mathbb{V}}
\newcommand{\bbW}{\mathbb{W}}
\newcommand{\bbX}{\mathbb{X}}
\newcommand{\bbY}{\mathbb{Y}}
\newcommand{\bbZ}{\mathbb{Z}}

\newcommand{\bbQp}{\mathbb{Q}_p}
\newcommand{\bbQpp}{\bar{\mathbb{Q}}_p}
\newcommand{\bbFp}{\mathbb{F}_p}
\newcommand{\bbFpp}{\bar{\mathbb{F}}_p}
\newcommand{\bbQQ}{\bar{\mathbb{Q}}}
\newcommand{\frakM}{\mathfrak{M}}
\newcommand{\frakp}{\mathfrak{p}}
\newcommand{\frakq}{\mathfrak{q}}
\newcommand{\fraka}{\mathfrak{a}}
\newcommand{\frakc}{\mathfrak{c}}

\newcommand{\frakX}{\mathfrak{X}}
\newcommand{\frakY}{\mathfrak{Y}}

\newcommand{\calA}{\mathcal{A}}
\newcommand{\calB}{\mathcal{B}}
\newcommand{\calC}{\mathcal{C}}
\newcommand{\calD}{\mathcal{D}}
\newcommand{\calE}{\mathcal{E}}
\newcommand{\calF}{\mathcal{F}}
\newcommand{\calG}{\mathcal{G}}
\newcommand{\calH}{\mathcal{H}}
\newcommand{\calI}{\mathcal{I}}
\newcommand{\calJ}{\mathcal{J}}
\newcommand{\calK}{\mathcal{K}}
\newcommand{\calL}{\mathcal{L}}

\newcommand{\calM}{\mathcal{M}}
\newcommand{\calN}{\mathcal{N}}

\newcommand{\calO}{\mathcal{O}}
\newcommand{\calP}{\mathcal{P}}
\newcommand{\calQ}{\mathcal{Q}}
\newcommand{\calR}{\mathcal{R}}
\newcommand{\calS}{\mathcal{S}}

\newcommand{\calT}{\mathcal{T}}
\newcommand{\calU}{\mathcal{U}}
\newcommand{\calV}{\mathcal{V}}
\newcommand{\calW}{\mathcal{W}}
\newcommand{\calX}{\mathcal{X}}
\newcommand{\calY}{\mathcal{Y}}
\newcommand{\calZ}{\mathcal{Z}}
\newcommand{\scrE}{\mathscr{E}}
\newcommand{\bbmi}{\mathbbm{i}}
\newcommand{\bbmj}{\mathbbm{j}}
\newcommand{\bbmk}{\mathbbm{k}}

\newtheorem{theorem}{Theorem}[section]
\newtheorem{lemma}[theorem]{Lemma}
\newtheorem{proposition}[theorem]{Proposition}
\newtheorem{definition}[theorem]{Definition}
\newtheorem{corollary}[theorem]{Corollary}
\newtheorem{conclusion}[theorem]{Conclusion}
\newtheorem{conjecture}[theorem]{Conjecture}
\newtheorem{remark}[theorem]{Remark}
\newtheorem{example}[theorem]{Example}
\newtheorem{notation}[theorem]{Notation}
\newtheorem{convention}[theorem]{Convention}
\newtheorem{problem}[theorem]{Problem}
\newtheorem{assumption}[theorem]{Assumption}
\newcommand{\bthm}{\begin{theorem}}
\newcommand{\ethm}{\end{theorem}}
\newcommand{\blem}{\begin{lemma}}
\newcommand{\elem}{\end{lemma}}
\newcommand{\bprop}{\begin{proposition}}
\newcommand{\eprop}{\end{proposition}}
\newcommand{\bdefn}{\begin{definition}}
\newcommand{\edefn}{\end{definition}}
\newcommand{\brmk}{\begin{remark}}
\newcommand{\ermk}{\end{remark}}
\newcommand{\bcor}{\begin{corollary}}
\newcommand{\ecor}{\end{corollary}}
\newcommand{\beg}{\begin{example}}
\newcommand{\eeg}{\end{example}}
\newcommand{\bitem}{\begin{itemize}}
\newcommand{\eitem}{\end{itemize}}
\newcommand{\rank}{{\rm rank}}
\newcommand{\isomap}{\stackrel{\thicksim}{\longrightarrow}}

%\thanks will become a 1st page footnote.
\thanks{This work is supported by the SFB/TR 45 `Periods, Moduli
Spaces and Arithmetic of Algebraic Varieties' of the DFG}

\begin{abstract}
This paper is devoted to the study of the Higgs bundle associated
with the universal abelian variety over the good reduction of a
Shimura curve of PEL type. Due to the endomorphism structure, the
Higgs bundle decomposes into the direct sum of Higgs subbundles
of rank two. They are basically divided into two types:
\emph{uniformizing type} and \emph{unitary type}. As the first
application we obtain the mass formula counting the number of
geometric points of the degeneracy locus in the Newton polygon
stratification. We show that each Higgs subbundle is Higgs
semistable. Furthermore, for each Higgs subbundle of unitary type,
either it is strongly semistable, or its Frobenius pull-back of a
suitable power achieves the upper bound of the instability. We
describe the Simpson-Ogus-Vologodsky correspondence for the Higgs
subbundles in terms of the classical Cartier descent.
\end{abstract}
%\fi

\maketitle

\section{Introduction}
Let $D$ be a quaternion division algebra over a totally real field
$F$ which is exactly split at one infinite place of $F$. By choosing
additionally a totally imaginary quadratic field extension $K$ of
$F$, the data $(D,K)$ allows one to define a Shimura curve of PEL
type (see \cite{De1}). In this paper, we study the Higgs bundle
$(E,\theta)$ associated with the universal abelian variety over
$\calM_0$, which is one of the geometrically connected components of
the good reduction of this Shimura curve modulo $p$. The passage of
the Higgs bundle from char 0 to char $p$ has two aims. The first is
to study the Newton polygon stratification of the moduli space in
char $p$. A similar method has already been extensively employed in
recent years (for example, see \cite{GO}, \cite{EvdG}, \cite{LSZ}, and so on).
The prototype of such study is the supersingular locus in the moduli
space of elliptic curves and the classical Deuring formula. From this
example one sees a basic phenomenon occurring in the geometry of
a moduli space in char $p$, namely the degeneration of the relative
Frobenius morphism along certain algebraic sublocus of the whole
moduli space. The second one is to investigate the relation between the
Higgs bundles over a char $p$ (or $p$-adic) field and the topology of the
underlying spaces in a char $p$ (or $p$-adic) field. In the classical
situation, that is, over the field of complex numbers, this is
beautifully expressed in the work of C. Simpson (see \cite{Si}).
Recently, analogous theories over a char $p$ or $p$-adic ground
field have emerged (see \cite{OV}, \cite{DW} and \cite{Fa}). We
intend to apply these new theories
to study the Higgs bundles over Shimura curves of PEL type in char $p$ and mixed characteristic. \\

Our results are built on the previous work on the Shimura curves
of PEL type, particularly the work of \cite{Car}, \cite{De1} and the
book \cite{Reim}. Let $\calM_0$ be the good
reduction in char $p$ of the Shimura curve of PEL type associated
with the quaternion division algebra $D$ and the imaginary quadratic
field $K$ (see \S2 for details), and let $f_0:\calX_0\to \calM_0$ be
the universal abelian variety and $(E,\theta)$ be the associated Higgs
bundle. Let $g$ be the genus of the Shimura curve $\calM_0$ which is
strictly greater than one by our choice of the level structure (see the
end of \S2). Most of the other notations appearing in the following
are collected at the end of this section.
\begin{theorem}[Proposition \ref{decomposition} and \ref{maximal Higgs field in char p}]\label{thm 1}
The Higgs bundle $(E,\theta)$ decomposes into the direct sum of rank two
Higgs subbundles:
$$
(E,\theta)=\bigoplus_{\phi\in \Phi}(E_{\phi},\theta_{\phi})\oplus
(E_{\bar \phi},\theta_{\bar \phi}),
$$
where the endomorphism subalgebra $\calO_{LK}\subset \calO_B$ acts
on the summand $E_{\phi}$ (resp. $E_{\bar \phi}$) via the character
$\phi \mod p$ (resp. $\bar \phi \mod p$). Assume further that $p\geq
2g$. Then for each $\phi\in \Phi$ (resp. $\bar \phi \in \bar \Phi$) with
$\phi|_F=\tau$ (resp. $\bar \phi|_F=\bar \tau$), the Higgs subbundle
$(E_{\phi},\theta_{\phi})$ (resp. $(E_{\bar \phi},\theta_{\bar
\phi})$) is of maximal Higgs field (see \cite{VZ}). Each of the remaining
Higgs subbundles in the above decomposition is of trivial (or equivalently, zero) Higgs
field.
\end{theorem}
The Higgs subbundles of maximal Higgs field are called of
\emph{uniformizing type}; while those of zero Higgs field are called of
\emph{unitary type}. In char 0, a Higgs bundle of uniformizing type
provides the uniformization of the base Shimura curve (see
\cite{VZ}) while that of unitary type corresponds to a unitary
representation of the topological fundamental group. We then analyze
the behavior of the iterated Frobenius morphism on the $(0,1)$-component
of each Higgs subbundle and derive the following results.

\begin{theorem} [Corollary \ref{slope} and Theorem \ref{mass formula}]\label{thm 2}
There are only two types of Newton polygons in $\calM_0(\bbF)$. Let
$\calS$ be the jumping locus of the Newton polygons. Then one has the
following formula in the Chow ring of $\calM_0$:
$$
\calS=\frac{1}{2}(1-p^{[F_{\mathfrak{p}}:\bbQ_p]})c_1(\calM_0).
$$
Taking the degree, one obtains the following mass formula for the Shimura
curve $\calM_0$:
$$
|\calS|=(1-p^{[F_{\mathfrak{p}}:\bbQ_p]})(1-g).
$$
\end{theorem}
From this formula one sees that the number of closed points in the
jumping locus of the Newton polygons is proportional to the topological
euler characteristic of the Shimura curve $\calM_0$. We would like to make the
following conjecture.

\begin{conjecture}
Let $\calM$ be the Shimura curve of Hodge type defined by
$G_{\bbQ}\to GSp_{\bbQ}$ with a suitable level structure, where
$G_{\bbQ}$ is the $\bbQ$-group of the units of a quaternion division
algebra $D$ over a totally real field $F$ by restriction of scalars.
Let $p$ be a prime number such that $\calM_0$ is the good reduction
of $\calM$ modulo $p$. Then there are exactly two possible Newton
polygons for the closed points of $\calM_0(\bbF)$ and the cardinality of the jumping
locus of the Newton polygons is equal to $\frac{1}{2}(1-p^d) \chi_{top}( \calM_0(\bbC))$,
where $d$ is the local degree $[F_{\mathfrak p}:\bbQ_p]$ and
$F_{\mathfrak p}$ is a splitting field of $D$ over $p$.
\end{conjecture}

In the case of Mumford's families of abelian varieties one sees much
evidence of the above conjecture by the work of R. Noot. He studied
the potential good reduction of a single abelian variety in a
Mumford's family, as well as its possible Newton polygons, using
Fontaine's theory (see \cite{No} and references therein). From the
proof of the mass formula one notices that certain Higgs subbundles
of unitary type have no contribution to the jumping of the Newton
polygons since the iterated relative Frobenius morphisms do not
degenerate on these subbundles. It turns out that these Higgs
subbundles of unitary type carry the extra property of being
strongly semistable, and this property characterizes the Higgs
subbundles of unitary type with non-degenerate iterated relative
Frobenius actions. This motivates us to study the (Higgs-)stability
of the Higgs subbundles under the Frobenius pull-backs in general.
For a semistable vector bundle $E$ over a smooth projective curve
$C$ in char $p$, the invariant $\nu(F_C^*E)$ (see \S\ref{stability}
for the definition) measures the extent of the instability under the
Frobenius pull-back. It is non-negative by definition and it is zero
if and only if $F_C^*E$ is still semistable. $E$ is strongly
semistable if $\nu(F_C^{n*}E)=0$ for all $n\geq 1$. It is well-known
that $\nu(F_C^*E)$ has the upper bound $(\rank(E)-1)(2g(C)-2)$ (see
Theorem \ref{instability}). Thus the extreme opposite of the
strongly semi-stability is the case that
$$
\nu(F_C^{n*}E)=0, \ 1\leq n\leq n_0,\ \textrm{and} \
\nu(F_C^{n_0*}E)=(\rank(E)-1)(2g(C)-2).
$$
We write the $Gal(L|\bbQ)$-orbit of $\Phi$ containing the
uniformizing place $\tau$ as follows:
$$
Hom_{\bbQp}(L_{\frakp},\bbQpp)=\{\phi_1,\cdots,\phi_d,\phi_1^*,\cdots,\phi_d^*\}
$$
with $\phi_1|_F=\phi_1^*|_F=\tau$ such that the Frobenius automorphism
$\sigma\in Gal(L_{\frakp}|\bbQ_p)$ acts on the orbit via the cyclic permutation. Thus we have the
following theorem.

\begin{theorem}[Proposition \ref{app-Ogus-V}, \ref{instable}, and \ref{sss}] \label{thm 3}
Assume that $p\geq 2g$. Then the following statements are true:
\begin{itemize}
  \item [(i)] Each Higgs subbundle in Theorem \ref{thm 1} is Higgs semistable of slope zero.
  In particular, the Higgs subbundles of unitary type are semistable.
  \item [(ii)] For $\phi\notin Hom_{\bbQp}(L_{\frakp},\bbQpp)$, the Higgs subbundle $(E_{\phi},\theta_{\phi})$
  of unitary type is strongly semistable (even \'{e}tale trivializable).
  The same is true for its bar counterpart.
  \item [(iii)] For $\phi_i \in Hom_{\bbQp}(L_{\frakp},\bbQpp)$ with $i\neq 1$, the Higgs subbundle
  $(E_{\phi_i},\theta_{\phi_i})$
  of unitary type satisfies:
  $$
\nu(F_{\calM_0}^{n*}E_{\phi_i})=0, \ 1\leq n\leq d-i, \mbox{ and }
\nu(F_{\calM_0}^{d-i+1*}E_{\phi_i})=2g-2.$$
The same is true for the $\phi_i^*$-summand with $i\neq 1$ and for
its bar counterpart.
\end{itemize}
\end{theorem}
For several reasons we are motivated to examine the
Simpson-Ogus-Vologodsky correspondence (see \cite{OV}) for the Higgs
subbundles in Theorem \ref{thm 1}. Let $(\calH^1_{dR},\nabla)$ be
the first relative de Rham bundle of $f_0$ with the canonical
Gauss-Manin connection, which has also an eigen-decomposition under the $\calO_{LK}$-action
(see Proposition \ref{decomposition}).
Let $F_{con}$ be the conjugate filtration on $\calH^1_{dR}$, which
is equally as important as the Hodge filtration in char $p$ geometry.

\begin{theorem}[Theorem \ref{Cartier descent} and Corollary \ref{HN filtration vs Hodge filtration}]\label{thm 4}
Assume that $p\geq \max\{2g,2([F:\bbQ]+1)\}$. Then for each $\phi\in \Phi$,
the Cartier transform of the direct summand
$(\calH^1_{dR,\phi},\nabla_{\phi})$ of $(\calH^1_{dR},\nabla)$ is
just the Cartier descent of
$Gr_{F_{con}}(\calH^1_{dR,\phi},\nabla_{\phi})$. The same is true
for its bar counterpart. As a consequence, for $\phi_i$ in Theorem
\ref{thm 3} (iii), the Harder-Narasimhan filtration on
$F_{\calM_0}^{d-i+1*}E_{\phi_i}$ is identified with the Hodge
filtration on $\calH^1_{dR,\phi_1^*}$. It is similar for the star and bar
counterparts.
\end{theorem}

By the above theorem one sees that the non-strongly semistable
Higgs subbundles of unitary type are closely related to the
Higgs subbundles of uniformizing type. In some sense one should
consider these two types of Higgs subbundles as the same one.
Compared with its char 0 analogue, the topological meaning of the
Higgs subbundles of uniformizing type is still unclear to us.

The paper is organized as follows. In \S\ref{reduction} the
construction of a Shimura curve of PEL type is briefly reviewed. In
\S\ref{local} some known results about Dieudonn\'{e} modules of the
abelian varieties corresponding to the points on $\calM_0(\bbF)$
are summarized. In \S\ref{global} the decomposition of the Higgs
bundle and the basic properties of the Higgs subbundles are
established. Applying the results in \S\ref{local} and
\S\ref{global}, we obtain the mass formula for the Shimura curve
$\calM_0$ in \S\ref{Newton jumping locus}. In \S\ref{stability} the
Higgs semi-stability as well as the semi-stability under Frobenius
pull-backs of the Higgs subbundles are discussed. The description of
the Simpson-Ogus-Vologodsky correspondence for the Higgs subbundles
is contained in \S\ref{SOV correspondence for Higgs subbundles}.\\

{\bf Notations and Conventions.}
\begin{itemize}
  \item [(i)] For a prime $\mathfrak{q}$ of a number field $E$, $E_{\mathfrak{q}}$
means the completion of $E$ with respect to $\mathfrak{q}$.  For a
field $E$ of char 0 (local or global), $\calO_E$ is the ring of
integers in $E$ and $\bar E$ is an algebraic closure of $E$.
$\bbQ_p^{ur}$ is the maximal unramified subextension of $\bbQ_p$.
Denote by $k$ a finite field of char $p$ and by $\bbF$ an
algebraic closure of $k$. Let $\sigma \in Gal(\bar \bbF|\bbF_p)$ be the
Frobenius automorphism, defined by $x\mapsto x^p$. It is restricted to
the Frobenius automorphism of $k$. For $\bbF$, denote by $W(\bbF)$ the
ring of Witt vectors and one has the canonical lifting of the
Frobenius automorphism of $\bbF$ to $W(\bbF)$, which is again denoted
by $\sigma$. It is similar for $W(k)$ and $k$.
  \item [(ii)] In this paper, $F$ is a fixed totally real number field of degree
$n\geq 2$, and $p$ is a rational prime number which is unramified in
$F$. $K$ and $L$ are two fixed imaginary quadratic field extensions
of $F$ (see \S2 for details). We put $\Psi = Hom_{\bbQ}(F,\bbR)$ and
$\Phi=Hom_{\bbQ}(L,\bar \bbQ)$. $D$ is a fixed quaternion division
algebra over $F$, which is exactly split at one infinite place $\tau\in
\Psi$ of $F$.
\iffalse
  \item [(iii)] For an abelian variety $A$ over $k$, $(\bbD(A),\calF, \calV)$ is the
associated Dieudonn\'{e} module, where $\calF$ (resp. $\calV$) is
the Frobenius morphism (resp. the Verschiebung morphism). The pair
$(\bbD(A),\calF)$ is in particular a F $\sigma$-crystal in the sense
of \cite{Ka1}.
\fi
  \item [(iii)] For an algebraic variety $X$ over $k$, one denotes by $F_X$ the
absolute Frobenius morphism. For a morphism $f: X\to Y$ over $k$ one
has the following commutative diagram of Frobenius morphisms:
$$
\xymatrix{ X\ar[dr]_{f}\ar[r]^{\calF_{X|Y}}&   X^{'} \ar[d]^{f^{'}}
\ar[r]^{\pi_{X|Y}}
                & X \ar[d]^{f}  \\
&  Y \ar[r]_{F_{Y}}
                & Y             }
$$
where the square in the diagram is the fiber product, $\calF_{X|Y}$ is the
relative Frobenius morphism and $\pi_{X|Y}\circ \calF_{X|Y}=F_X$.
For a vector bundle $\calE$ over $X$ in char $p$, sometimes we denote
the (iterated) Frobenius pull-back $F_X^{*n}\calE (n\geq 1)$ by $\calE^{(p^{n})}$.
  \item [(iv)] In this paper, the term `\emph{reduction modulo $p$}' means the
following: let $R$ be a DVR of mixed characteristic $(0,p)$ with the
residue field $k(R)$, and $M$ be an object defined over $R$, which can
be a module or a scheme. Then the reduction of $M$ modulo $p$ is the
base change of $M$ from $R$ to $k(R)$.
\end{itemize}

\textbf{Acknowledgements.} We would like to thank the referee for his/her careful reading of our paper
and helpful advice. We thank C. Deninger for useful discussions on \cite{Den} and \cite{DW}.
%The Higgs subbundles in Theorem \ref{thm 3} (ii) are subject to his construction with A. Werner in \cite{DW}.
%The resulting $p$-adic
%representations of geometric \'{e}tale fundamental group of the
%generic fiber are actually isomorphic to the original ones from
%geometry. We would like to clarify this point together with certain
%generalizations to a higher dimensional base in a forthcoming paper.
Special thanks go to A. Langer for his useful comments on \S\ref{stability} and particularly
the clarification of a main result in \cite{LS} (see Theorem \ref{result of LS}). \\

\section{Quaternion division algebras and the good reduction of a Shimura curve}
\label{reduction}

Let $D$ be a quaternion division algebra over $F$, which is split at
the infinite place $\tau$ and ramified at all remaining infinite places.
That is, one has the following isomorphisms:
$$D\otimes_{F,\tau}\bbR\cong M_2(\bbR), \mbox{ and } D\otimes_{F,\psi}\bbR\cong
\bbH \mbox{ for } \psi\neq \tau,$$ where $\bbH$ is the Hamiltonian quaternion
algebra over $\bbR$. One considers the $F$-group of the units
$D^{\times}$ as a $\bbQ$-group by restriction of scalars and defines a
homomorphism of real algebraic groups
\begin{eqnarray*}
h_D:\mathbb{S}=\mathbb{C}^{\times}&\rightarrow&
D^{\times}(\bbR)\cong GL_2(\bbR)\times (\bbH^{\times})^{n-1}\\
 z=x+iy&\mapsto&(\begin{pmatrix}x & y\\ -y&x\end{pmatrix},1,\cdots,1).
\end{eqnarray*}
The $D^{\times}(\bbR)$-conjugacy class $X$ of $h_D$ defines a
Shimura curve $Sh_D$ over the reflex field $F$, where $F$ is considered as a
subfield of $\bbC$ via the embedding $\tau$. For every open compact
subgroup $C\subset D^{\times}(\mathbb{A}^f)$, $Sh_{D,C}=Sh_D/C$ is a
projective curve over $F$, and one has the identification of its
complex points
$$
Sh_{D,C}(\mathbb{C})=D^{\times}(\bbQ)\setminus(X\times
D^{\times}(\mathbb{A}^f))/C
$$
where $\bbA_f$ is the ring of finite ad\`{e}les of $\bbQ$ and
$D^{\times}$ acts on $X$ by the conjugation and on
the second summand by the left multiplication. \\

In certain cases $Sh_{D,C}(\bbC)$ is known to parameterize the
principally polarized abelian varieties over $\bbC$ with special
Mumford-Tate groups (see for example \cite{Mum} \S4 and \cite{VZ} \S5).
It belongs to the category of Shimura varieties of Hodge type. In
this paper we are going to study a related Shimura curve $Sh_{G}$
which is of PEL type. However the Shimura curve of Hodge type
provides the motivation for the further study of the current paper.
Now we recall the so-called `mod\`{e}le \'{e}trange'
construction in \cite{De1}. First we need to choose an imaginary
quadratic field $\bbQ(\alpha)$ with $\alpha\in \bbC$ such that $p$
is split in it. We put the composite $K=F(\alpha)$ and it will be
fixed in the whole discussion. Considering $F^{\times}$ and
$K^{\times}$ as $\bbQ$-groups, we define a new $\bbQ$-group $G$ by
the following short exact sequence
$$
1\rightarrow F^{\times}\rightarrow D^{\times}\times
K^{\times}\stackrel{\pi}{\rightarrow} G\rightarrow 1,
$$
where $F^{\times}\rightarrow D^{\times}\times K^{\times}$ is given
by $f\mapsto (f,f^{-1})$. We fix a subset $\Psi_K\subset
Hom_{\bbQ}(K,\bbC)$ which induces a bijection to $\Psi$ by
restriction to $F$. Note that $\Psi_K$ is obtained by the trivial extensions of
all embeddings of $F$ into $\bbC$ to embeddings of $K=F+F\alpha$ into $\bbC$. One has an
identification $K^{\times}(\bbR)\cong \prod_{\psi\in \Psi}
\bbC^{\times}$ and defines
$$
h_K:\bbS\rightarrow K^{\times}(\bbR)=\bbC^{\times}\times
\prod_{\psi\neq \tau}\bbC^{\times},\ z\mapsto (1,z,\cdots,z).
$$

Let $X'$ denote the conjugacy class of
$$
h_{G}=\pi_{\mathbb{R}}\circ(h_D\times h_K):\mathbb{S}\rightarrow
G(\mathbb{R}).
$$
It defines a Shimura curve $Sh_{G}$ over $K$, where $K$ is a
subfield of $\bbC$ via the map $\tau\in \Psi\cong \Psi_K$. A compact
open subgroup $C$ of $G(\bbA_f)$ defines a projective curve
$Sh_{G,C}$. For a suitable $C'\subset D^{\times}(\bbA_f)$, the neutral
component of $Sh_{D,C'}$ and $Sh_{G,C}$ are isomorphic to each other
over certain number field (see \S1 of \cite{Reim}). The Shimura
curve $Sh_{G,C}$ parameterizes the isogeny classes of abelian
varieties over $K$ with PEL structure which we describe briefly as
follows. Let $B=D\otimes_F K$. Define the natural involution on $B$
by the formula
$$
(x\otimes y)'=x^*\otimes \bar{y},
$$
where $*$ is the main involution on $D$ and $^{-}$ is the complex
conjugation on $K$ which is the generator of $Gal(K|F)$. Let $V$ be
the underlying $\bbQ$-vector space of $B$. There exists a
non-degenerate alternating $\bbQ$-bilinear form
$$
\Theta: V\times V\rightarrow \bbQ$$ such that for all $b\in B,x,y\in
V$, $\Theta(bx,y)=\Theta(x,b'y)$. It turns out that $G$ is the group
of $B$-module automorphisms of $V$
preserving the bilinear form. So one has the natural linear representation
$\xi_{\bbQ}:G(\bbQ)=Aut_{B}(V,\Theta)\subset Aut_{\bbQ}(V)$. \\

We write $p\calO_F=\prod_{i=1}^{r}\mathfrak{p}_i$. By fixing an
embedding $\bar \bbQ\to \bar \bbQ_p$, one obtains an bijection
between $\Psi=Hom_{\bbQ}(F,\bar \bbQ)$ and
$\coprod_{i=1}^{r}Hom_{\bbQ_p}(F_{\mathfrak{p}_i},\bar \bbQ_p)$.
After a rearrangement of indices we can assume that, under the above
bijection, $\tau$ lies in $Hom_{\bbQ_p}(F_{\mathfrak{p}_1},\bar
\bbQ_p)$. We fix the notation $\mathfrak{p}=\mathfrak{p}_1$ for the
whole paper. Since $p$ is split in $\bbQ(\alpha)$ by assumption,
$\mathfrak{p}_i\calO_K=\mathfrak{q}_i \bar{\mathfrak{q}}_i$ for
each $i$, where the two primes of $K$ over $\mathfrak p_i$ are
distinguished in such a manner that $\Psi_K$ is bijectively mapped onto
$\coprod_{i=1}^{r}Hom_{\bbQ_p}(K_{\mathfrak{q}_i},\bar \bbQ_p)$
under the previous identification map. Now we fix a totally imaginary
quadratic extension $L$ of $F$ which is contained in $D$. Then $L$
splits $D$ globally. Namely there is an isomorphism of $L$-algebras
$D\otimes_F L\cong M_2(L)$. Furthermore one can assume that each
$\mathfrak{p}_i$ stays prime in $L$. So the composite field $LK$ is
particularly unramified over $p$. One writes the prime ideal
decomposition as
$$
p\calO_{LK}=\prod_{i=1}^{r}\mathfrak{q}_i \bar{\mathfrak{q}}_i,
$$
and one has a natural isomorphism of $\bbQ_p$-algebras
$$
LK\otimes_{\bbQ}\bbQ_p\cong
\prod_{i=1}^{r}LK_{\mathfrak{q}_i}\times
LK_{\bar{\mathfrak{q}}_i}.
$$
Moreover for each $i$ one has an isomorphism
$LK_{\mathfrak{q}_i}\otimes_{\bbQ_p}\bbQ_p^{ur}\cong
\prod_{Hom_{\bbQ_p}(LK_{\mathfrak{q}_i},\bar \bbQ_p)}\bbQ_p^{ur}$.
It is similar for the bar counterpart. Then we obtain an isomorphism of
$\bbQ_p^{ur}$-algebras
$$
LK\otimes_{\bbQ}\bbQ_p^{ur}\cong
\prod_{i=1}^{r}(\prod_{Hom_{\bbQ_p}(LK_{\mathfrak{q}_i},\bar
\bbQ_p)}\bbQ_p^{ur}\times
\prod_{Hom_{\bbQ_p}(LK_{\bar{\mathfrak{q}}_i} ,\bar
\bbQ_p)}\bbQ_p^{ur} ).
$$
It induces on the rings of integers a $\bbZ_p$-algebra isomorphism.
One can simplify the notations by using the identification $$
Hom_{\bbQ}(LK,\bar
\bbQ)=\coprod_{i=1}^{r}Hom_{\bbQ_p}(LK_{\mathfrak{q}_i},\bar
\bbQ_p)\times Hom_{\bbQ_p}(LK_{\bar{\mathfrak{q}}_i},\bar
\bbQ_p),
$$
and the partition $Hom_{\bbQ}(LK,\bar \bbQ)=\Phi\coprod \bar \Phi$,
where $\Phi=Hom_{\bbQ}(L,\bar \bbQ)$ is identified with the subset
of $Hom_{\bbQ}(LK,\bar \bbQ)$ by extending each embedding of $L$ into
$\bar \bbQ$ to an embedding of $LK=L(\alpha)=L+L\alpha$ into $\bar
\bbQ$ which is the identity on $\alpha$. Thus we write the above isomorphism
of $\bbZ_p$-algebras in the form:
\begin{eqnarray}\label{action-char}
\prod_{\phi\in\Phi}w(\phi)\times \bar{w}(\phi):\calO_{L
K}\otimes_{\mathbb{Z}}W(\bbF) \isomap \prod_{\phi\in\Phi}
W(\bbF)\times W(\bbF),
\end{eqnarray}
where for each $\phi\in \Phi$, $w(\phi)\in
Hom_{\bbQ_p}(LK_{\mathfrak{q}_i},\bar \bbQ_p)$ and $\bar{w}(\phi)
\in Hom_{\bbQ_p}(LK_{\bar{\mathfrak{q}}_i},\bar \bbQ_p)$ for
certain $i$. By abuse of notations we also write the character
$w(\phi)$ (resp. $\bar{w}(\phi)$) as $\phi$ (resp. $\bar
\phi$) simply. In the following we come to an important notion for this section.

\begin{definition}\label{sheaf of certain type}
Let $S$ be an $\calO_{F_\mathfrak{p}}$-scheme and $\calE$ be a
locally free coherent $\calO_{S}$-sheaf. It is said to be a sheaf of
type $(L,\Psi_K)$ if $\calO_{LK}\subset End_{\calO_S}(\calE)$ and
$\calE\otimes_{\bbZ} W(\bbF)$ has a decomposition induced by the
isomorphism (\ref{action-char}) as follows:
\begin{eqnarray}\label{sheaf-decomp}
\calE\otimes W(\bbF)=\bigoplus_{\phi\in \Phi} (\calE_{\phi}\oplus
\calE_{\bar\phi}),
\end{eqnarray}
where $\calE_{\phi}$ corresponds to the character $w(\phi)$ and
$\calE_{\bar\phi}$ corresponds to the character $\bar{w}(\phi)$ with the rank
condition: $\calE_{\phi}$ and $\calE_{\bar\phi}$ are of rank one if
$\phi|_F=\tau$, while $\calE_{\phi}$ is of rank two and $\calE_{\bar\phi}=0$
if $\phi|_F\neq \tau$.
\end{definition}

In order to define a level structure one shall choose an integral
structure of the $\bbQ$-vector space $V$. One chooses an order $\calO_D$
of $D$ containing $\calO_L$ with certain additional properties (see
\S2 in \cite{Reim}), and put $\calO_B=\calO_D\otimes_{\calO_F}
\calO_K$ (so $\calO_{LK}\subset \calO_B$). Then one takes the lattice
$V_{\bbZ}$ of $V$ to be the free $\bbZ$-module $\calO_B$ and puts
$G(\bbZ)=Aut_{\calO_B}(V_{\bbZ}, \Theta)$. Thus one has an
integral structure $\xi: G(\bbZ)\to Aut_{\bbZ}(V_{\bbZ})$ of the
$\bbQ$-algebraic group morphism $\xi_{\bbQ}$.

\begin{proposition}[Proposition 2.14 and Corollary 3.14 in \cite{Reim}]\label{local-ring}
For every level structure $C=C_p\times C^p\subset G(\bbA_f)$ with
$C_p=G(\bbZ_p)$ and $C^p$ small enough, there exists a proper
$\calO_{F_{\mathfrak{p}}}$-scheme $\calM_{C}$ which is the coarse
moduli space of certain moduli functor of PEL type (see Proposition
2.14 in \cite{Reim} for the description of the moduli functor) with
the endomorphism algebra $\calO_B$. Furthermore, if $D$ is assumed
to be split at $\mathfrak{p}$, then the reduction $\calM_C$ modulo
$p$ is smooth over $\bbF$.
\end{proposition}

We take one of the geometrically connected components $\calM$ of
$\calM_{C}$ with the reduction $\calM_0$ modulo $p$. For our purpose
we shall take $C^p$ small enough so that we have the
universal abelian scheme $f:\calX\to \calM$. Under this assumption
the genus of $\calM$ must be strictly greater than one. By the
construction of the moduli functor, the injection
$\calO_B\hookrightarrow End_{\calM}(\calX)$ turns
$R^1f_{*}\calO_{\calX}$ into a sheaf of $(L,\Psi_K)$-type.

\section{Dieudonn\'{e} modules and Newton polygons}\label{local}

Let $A$ be an abelian variety which is represented by an
$\bbF$-rational point of $\calM_0$. Let $(\bbD=\bbD(A),\calF,\calV)$
be the associated (contravariant) Dieudonn\'{e} module. $\bbD$ is a
free $W(\bbF)$-module of rank $8n$ and one has the identifications
of $k$-vector spaces:
$$
\bbD/p\bbD=H_{dR}^1(A),\ \calV\bbD/p\bbD=H^0(A,\Omega_A),\mbox{ and }
\bbD/\calV\bbD(A)=H^1(A,\calO_A).
$$
In this section we shall analyze the structure of $\bbD$ in the
presence of the endomorphism structure. Actually, since
$\calO_B\subset End(A)$, it follows that $\calO_B \subset End(\bbD)$
and particularly $\calO_{LK}\subset End(\bbD)$. Therefore $\bbD$ is
an $\calO_{LK}\otimes_{\bbZ}W(\bbF)$-module. The isomorphism
(\ref{action-char}) in \S2 gives the decomposition
$$
\bbD=\bigoplus_{\phi\in \Phi}(\bbD_{\phi}\oplus \bbD_{\bar{\phi}}).
$$
We put for $1\leq i\leq r$ the local degree
$f_i=[F_{\frakp_i}:\bbQp]$ and
$L_{\frakp_i}=L\otimes_{F}F_{\frakp_i}$. For an element
$\phi\in\Phi$ one defines $\phi^*\in \Phi$ to be the another element
whose restriction to $F$ is the same as that of $\phi$. The
following proposition contains the basic properties of each direct
summand in the above decomposition.
\begin{proposition} \label{Dieudonne module structure}
The Dieudonn\'{e} module $\bbD$ has the following properties:
\begin{itemize}
    \item [(i)] $\calO_{LK}$ acts on $\bbD_{\phi}$ (resp.
$\bbD_{\bar{\phi}}$) via the character $w(\phi)$ (resp.
$\bar{w}(\phi)$).
    \item [(ii)] There is an endomorphism $\Pi\in \calO_B\otimes
    \bbZ_p$ which induces a morphism $\Pi: \bbD_{\phi}\to
    \bbD_{\phi^*}$. It is an isomorphism for $\phi|_F=\tau$.
    \item [(iii)] For $\phi|_{F}\neq\tau$, $\calF(\bbD_{\phi})=
    \bbD_{\sigma\phi}$. For $\phi|_{F}=\tau$, $p\bbD_{\sigma\phi}\subsetneqq \calF \bbD_{\phi}\subsetneqq
\bbD_{\sigma\phi}$.
    \item [(iv)] For each $\phi\in \Phi$,
both $\bbD_{\phi}$ and $\bbD_{\bar{\phi}}$ are of rank $2$.
 \item
[(v)] The polarization induces a perfect alternative pairing between
$\bbD_{\phi}$ and $\bbD_{\bar{\phi}}$. Moreover $\bbD_{\phi}\perp
\bbD_{\phi'}$ unless $\phi'=\bar{\phi}$. Thus $\bbD_{\phi}$ and
$\bbD_{\bar{\phi}}$ are dual to each other.
\end{itemize}
\end{proposition}
\begin{proof}
(i) follows from the definition. The existence of $\Pi\in
\calO_B\otimes \bbZ_p$ with the property as in (ii) is actually a part
of the conditions on $\calO_B$ (see \S2 in \cite{Reim}). Clearly
$\calF$ commutes with the $\calO_B$-action on $\bbD$. Since $\calF$ is
$\sigma$-semilinear, one has $\calF(\bbD_{\phi})\subseteq
\bbD_{\sigma\phi}$ by (i). Similarly one has $\calV
(\bbD_{\sigma\phi}) \subseteq \bbD_{\phi}$. For a fixed $\phi$ we
consider the short exact sequence
$$
0\rightarrow
\frac{\calV(\bbD_{\sigma\phi})}{p\bbD_{\phi}}\rightarrow
\frac{\bbD_{\phi}}{p\bbD_{\phi}}\rightarrow
\frac{\bbD_{\phi}}{\calV\bbD_{\sigma\phi}} \rightarrow 0.
$$
By the rank condition in Definition \ref{sheaf of certain type}, one
has $\dim_{\bbF}\frac{\bbD_{\phi}}{\calV\bbD_{\sigma\phi}}$ is equal
to one if $\phi|_F=\tau$, and equal to two if $\phi|_F\neq \tau$. Moreover
$\dim_{\bbF}\frac{\bbD_{\bar \sigma\phi}}{\calV\bbD_{\bar \phi}}$ is
equal to one in the former case and zero in the latter case. By
duality, namely $H^0(A,\Omega_A)^*\cong H^1(A^t,\calO_{A^t})$ with
$A^t$ the dual abelian variety of $A$, it follows that
$\dim_{\bbF}\frac{\calV(\bbD_{\sigma\phi})}{p\bbD_{\phi}}$ is equal
to one in the former case and equal to zero in the latter case. So
in both cases $\dim_{\bbF}\frac{\bbD_{\phi}}{p\bbD_{\phi}}=2$ and therefore
$\rank_{W(\bbF)}\bbD_{\phi}=2$. (iv) follows from (v). By the above proof,
we have $\calV (\bbD_{\sigma\phi})=p \bbD_{\phi}$ for
$\phi|_F\neq\tau$ and $p\bbD_{\phi}\varsubsetneq
\calV(\bbD_{\sigma\phi})\varsubsetneq \bbD_{\phi}$ otherwise. By
applying $\calF$ to both sides and dividing by $p$ if necessary, one
obtains (iii). Finally, since $\psi(lx,y)=\psi(x,l'y)$ for all
$x,y\in \bbD$ and $l\in B$. We take two idempotents
$l_{\phi},l_{\phi'}$ with $l_{\phi}\in (\calO_{LK})_{\phi}$ and
$l_{\phi'}\in (\calO_{LK})_{\phi'}$. Then
$\psi(l_{\phi}x,l_{\phi'}y)=\psi(x,l_{\phi}'l_{\phi'}y)=0$ unless
$\bar{\phi}=\phi'$, since $l_{\phi}'\in (\calO_{LK})_{\bar{\phi}}$.
Thus (v) follows.
\end{proof}

In the following we determine the possible Newton polygons of
$\bbD$. As it is an isogeny invariant, we introduce the (F-)$\sigma$-isocrystal
$(N=\bbD\otimes_{\bbZ}\bbQ,\calF)$, and similarly for $\phi\in
\Phi$, the direct summand $N_{\phi}$ (resp. $N_{\bar \phi}$) which
itself is generally not a sub F-$\sigma$-isocrystal by the previous
proposition. For each $1\leq i\leq r$, we put $ N_i=\oplus_{\phi\in
Hom_{\bbQ_p}(L_{\frakp_i},\bar \bbQ_p)}N_{\phi} $ and similarly for
$N_{\bar i}$. Then $N_i$ and $N_{\bar i}$ are indeed F-$\sigma$-isocrystals for each $i$. However if we set
$\calF_i=(\calF^{f_i}\Pi)|_{N_i}$, then by the relation
$\sigma^{f_i}\phi=\phi^*$, one sees that $(N_{\phi},\calF_i)$ is
indeed an F-$\sigma^{f_i}$-isocrystal for each $\phi\in
Hom_{\bbQ_p}(L_{\frakp_i},\bar \bbQ_p)$, and similarly for the bar
counterpart.
\begin{proposition}
Let $(N_i, \calF)$ be the F-$\sigma$-isocrystal as above. Then it
has the following possible Newton slopes:
\begin{itemize}
    \item [(i)] For $i=1$ the Newton slopes are either $4f_1\times 1/2f_1$ or $2f_1\times
    (0,1/f_1)$.
    \item [(ii)] For $i\geq 2$ the Newton slope is $4f_i\times 0$.
\end{itemize}
\end{proposition}
\begin{proof}
Since $\calF: N_{\phi}\rightarrow N_{\sigma\phi}$ is an isogeny of
isocrystals, $(N_{\phi},\calF_i)$ and $(N_{\sigma\phi},\calF_i)$
have the same Newton slopes. So the computation is reduced to the
$\sigma^{f_i}$-isocrystals $(N_{\phi},\calF_i)$ of height 2.
Thus the result follows easily from the classification of
the isocrystals of height 2 over $\bbF$ (cf. Lemma 4.4 in \cite{Reim} or
\cite{Zink}).
\end{proof}

We put $d=[F_{\frakp}:\bbQ_p]=f_1$. The following corollary follows
easily from the last proposition.
\begin{corollary}\label{slope}
Let $A$ be an abelian variety which is represented by an
$\bbF$-rational point of $\calM_0$. Then the Newton polygon of $A$ is of the following two possible
types:
$$
(4n-2d)\times 0,\ 2d\times 1/d,\ 2d\times(1-1/d),\ (4n-2d)\times 1;
$$
and
$$
4(n-d)\times 0,\ 4d\times 1/2d,\ 4d\times (1-1/2d),\ 4(n-d)\times 1.
$$
\end{corollary}
\begin{proof}
It suffices to notice that for each $i$, $N_i$ and $N_{\bar{i}}$ are
dual to each other as $\sigma$-isocrystals by Proposition
\ref{Dieudonne module structure} (v).
\end{proof}
\begin{remark}
We see that there are only two possible Newton polygons for closed
$\bbF$-points of the moduli space. The existence of the abelian
varieties with the given Newton polygons was shown by Honda-Tate
theory. We refer to \cite{Car} or \cite{Reim} for the details.
\end{remark}

\section{The Decomposition of the Higgs bundle over a Shimura curve in char $p$} \label{global}
Let $f: \calX\to \calM$ be the universal abelian scheme in \S2. By
abuse of notations we denote it again by $f$ the base change to $\bar
\bbZ_p$. Let $f^0: \calX^0\to \calM^0$ be the base change of $f$ to
$\bar \bbQ_p$ and $f_0: \calX_0\to \calM_0$ be the base change to
$\bbF$. Let
$\calH^1_{dR}=R^1f_{0*}(\Omega^{\cdot}_{\calX_0|\calM_0},d)$ be the
first relative de Rham bundle over $\calM_0$. We put the first Hodge
bundle $E^{1,0}=f_{0*}\Omega^1_{\calX_0|\calM_0}$ and the second
Hodge bundle $E^{0,1}=R^1f_{0*}\calO_{\calX_0}$. By the
$E_1$-degeneration of the Hodge to de Rham spectral sequence, one
has the short exact sequence
$$
0\rightarrow E^{1,0} \rightarrow \calH^1_{dR}\rightarrow E^{0,1}
\rightarrow 0.
$$
It is well-known that $\calH^1_{dR}$ is endowed with the Gauss-Manin
connection $\nabla$. By taking the grading $(\calH^1_{dR},\nabla)$
with respect to the Hodge filtration,
we obtain the Higgs bundle in char $p$: $(E,\theta)=(E^{1,0}\oplus E^{0,1},
\theta^{1,0}\oplus \theta^{0,1})$ with $\theta^{0,1}=0$. \\

By construction, the endomorphism ring of the universal abelian
variety $\calX_0$ over $\calM_0$ contains $\calO_B$. Thus each
element $b\in \calO_B$ induces a morphism $b: \calX_0\to \calX_0$
over $\calM_0$. Let $\calO_{LK}\subset \calO_B$ be the maximal
abelian subgroup as in \S2. We have the following decomposition
under the $\calO_{LK}$-action.
\begin{proposition}\label{decomposition}
The first relative de Rham bundle $(\calH^1_{dR},\nabla)$ with the
Gauss-Manin connection admits a decomposition into the direct sum of
rank two subbundles with an integrable connection
$$
(\calH^1_{dR},\nabla)=\bigoplus_{\phi\in
\Phi}(\calH^1_{dR,\phi},\nabla_{\phi})\oplus (\calH^1_{dR,\bar
\phi},\nabla_{\bar \phi}),
$$
such that $\calO_{LK}$ acts on $\calH^1_{dR,\phi}$ (resp.
$\calH^1_{dR,\bar \phi}$) via the character $\phi \mod p$ (resp.
$\bar \phi \mod p$). It induces the decomposition of the Higgs
bundle into the direct sum of rank two Higgs subbundles
$$
(E,\theta)=\bigoplus_{\phi\in \Phi}(E_{\phi},\theta_{\phi})\oplus
(E_{\bar \phi},\theta_{\bar \phi}).
$$
Furthermore, by writing
$$
E_{\phi}=E^{1,0}_{\phi}\oplus E^{0,1}_{\phi} \mbox{ and } E_{\bar
\phi}=E^{1,0}_{\bar \phi}\oplus E^{0,1}_{\bar \phi},
$$
one has for $\phi|_F=\tau$, $\rank\ E^{0,1}_{\phi}=\rank\
E^{0,1}_{\bar \phi}=1$; while for $\phi|_F\neq \tau$, $\rank\
E^{0,1}_{\phi}=2$ and $\rank\ E^{0,1}_{\bar \phi}=0$.
\end{proposition}
\begin{proof}
The decomposition of $\calH^1_{dR}$ with respect to the
$\calO_{LK}$-action follows from Proposition \ref{Dieudonne module
structure} (i). Because $\calO_{LK}$ acts on $\calX_0$ as
endomorphisms over $\calM_0$, it induces an action on the relative
de Rham complex as endomorphisms of complexes. Taking the
hypercohomology, it induces an action on the Hodge filtration
$0\subset E^{1,0}\subset \calH^1_{dR}$. In other words, $E^{1,0}$ is
an $\calO_{LK}$-invariant subbundle of $\calH^1_{dR}$. Thus one has
the corresponding decomposition on $E^{1,0}$. $E^{0,1}$ is the
quotient bundle $\frac{\calH^1_{dR}}{E^{1,0}}$. Then for each $\phi$
(resp. $\bar \phi$), one has an injective morphism
$\frac{\calH^1_{dR,\phi}}{E^{1,0}_{\phi}}\to E^{0,1}$ (resp. for
$\bar \phi$) induced by $\calH^1_{dR} \twoheadrightarrow E^{0,1}$.
So one has an isomorphism
$$
E^{0,1}\cong \bigoplus_{\phi \in \Phi}
\frac{\calH^1_{dR,\phi}}{E^{1,0}_{\phi}}\oplus
\frac{\calH^1_{dR,\bar \phi}}{E^{1,0}_{\bar \phi}}.
$$
Denote $ \frac{\calH^1_{dR,\phi}}{E^{1,0}_{\phi}}$ by
$E^{1,0}_{\phi}$ (similarly for $\bar \phi$), we obtain the
decomposition of $E^{0,1}$. By the short exact sequence
$$
0\rightarrow {E^{1,0}_{\phi}}\rightarrow
\calH^1_{dR,\phi}\rightarrow{E^{0,1}_{\phi}}\rightarrow 0,
$$
the bundle $E_{\phi}=E^{1,0}_{\phi}\oplus E^{0,1}_{\phi}$ has the
same rank as $\calH^1_{dR,\phi}$, which is two by Proposition
\ref{Dieudonne module structure} (iv). It is similar for $\bar \phi$. It
is clear that the resulting decomposition on $E^{0,1}$ coincides
with the induced action of $\calO_{LK}$ on
$R^1f_{0*}\calO_{\calX_0}=E^{0,1}$ by taking the higher direct image.
Since the bundle $E^{0,1}$ is the modulo $p$ reduction of
$R^1f_{*}\calO_{\calX}$, it is a sheaf of $(L,\Psi_K)$-type. The
assertions about the ranks of $E^{0,1}_{\phi}$ and
$E^{0,1}_{\bar\phi}$ follow from the rank condition in Definition
\ref{sheaf of certain type}. Finally the $\calO_{LK}$-action
decomposes the Gauss-Manin connection as well. In fact, in char $0$
one can show that the $\calO_{LK}$-action on the relative de Rham
bundle is flat with respect to (in other words, commutes with) the Gauss-Manin connection because
the endomorphism algebra defines the flat Hodge cycles on the
relative Betti cohomology. By reduction modulo $p$, the
$\calO_{LK}$-action also commutes with $\nabla$. Because
$\calO_{LK}$ acts on the direct summands via the characters, $\nabla$
preserves each direct summand in the decomposition. The Higgs field
$\theta$ on $E$ decomposes accordingly.
\end{proof}

\begin{corollary}\label{decomposition in char 0}
The Hodge-to-de Rham spectral sequence of the relative de Rham
bundle \\ $R^1f_*(\Omega_{\calX|\calM},d)$ degenerates at $E_1$-level. By
taking the grading of $(R^1f_*(\Omega_{\calX|\calM},d),\nabla)$ with
respect to the Hodge filtration, one obtains the Higgs bundle $(\tilde
E,\tilde \theta)$ over $\calM$. The $\calO_{LK}$-action on the
universal abelian scheme $\calX$ over $\calM$ as endomorphisms
induces a decomposition of Higgs bundles
$$
(\tilde E,\tilde \theta)=\bigoplus_{\phi\in \Phi}(\tilde
E_{\phi},\tilde \theta_{\phi})\oplus (\tilde E_{\bar \phi},\tilde
\theta_{\bar \phi}).
$$
The modulo $p$ reduction of the above decomposition is the one over $\calM_0$ given
in Proposition \ref{decomposition}.
\end{corollary}
\begin{proof}
It suffices to show the $E_1$-degeneration of the Hodge-to-de Rham
spectral sequence of $f$. It is equivalent to show that the natural
morphism $f_*\Omega_{\calX|\calM}\to R^1f_*(\Omega_{\calX|\calM},d)$
is injective. By tensoring with $\bbQ$, the above morphism is
injective by the well-known $E_1$-degeneration of the Hodge-to-de Rham spectral
sequence for the first relative de Rham bundle of $f^0$ (the generic
fiber of $f$). So the kernel of the morphism consists of only
$p$-torsions. By modulo $p$ and the $E_1$-degeneration of the closed
fiber $f_0$ of $f$, there is actually no $p$-torsions. Thus the
Hodge-to-de Rham spectral sequence of $f$ degenerates at $E_1$-level as
well.
\end{proof}

Next we proceed to deduce some basic properties of the Higgs
subbundles from the above proposition. According to this proposition, the Higgs
subbundle $(E_{\phi},\theta_{\phi})$ (resp. $(E_{\bar
\phi},\theta_{\bar \phi})$) has two nontrivial parts, namely, the $(1,0)$-part and the
$(0,1)$-part, if and only if $\phi|_F=\tau$ (resp. $\bar \phi|_F=\bar \tau$).
It is clear that there are totally four such direct summands in the
decomposition. We consider them first.
\begin{proposition}\label{uniformizing Higgs bundle}
Let $\phi,\phi^*$ be two unique elements of $\Phi$ whose restriction
to $F$ is equal to $\tau$. Then one has an isomorphism of Higgs bundles
$(E_{\phi},\theta_{\phi})\cong (E_{\phi^*},\theta_{\phi^*})$. One
has also an isomorphism for the bar counterpart.
\end{proposition}
\begin{proof}
In case of $\phi|_F=\tau$ the endomorphism $\Pi\in \calO_B\otimes
\bbZ_p$ of $\calX_0$ over $\calM_0$ induces the endomorphism $\Pi
\in End(\calH^1_{dR})$ which is in fact an automorphism. By
restricting $\Pi$ to each closed point in $\calM_0$, one knows
from Proposition \ref{Dieudonne module structure} (ii) (modulo
$p$) that it induces an isomorphism $\Pi:
\calH^1_{dR,\phi}\rightarrow \calH^1_{dR,\phi^*}$. Since $\Pi$
commutes with the Gauss-Manin connection and the Hodge filtration,
it induces an isomorphism of Higgs bundles by taking the grading with
respect to the Hodge filtration:
$$
\Pi: (E_{\phi},\theta_{\phi})\cong (E_{\phi^*},\theta_{\phi^*}).
$$
\end{proof}
The following result asserts that the chern class of the base
Shimura curve $\calM_0$ is in fact represented by the second Hodge
bundle of the Higgs subbundles appearing in the above proposition.
This is one of significant features of the above Higgs subbundles.

\begin{proposition}\label{maximal Higgs field in char p}
Assume that $p\geq 2g$. Then for $\phi\in \Phi$ with $\phi|_{F}=\tau$,
the Higgs bundle $(E_{\phi},\theta_{\phi})$ in char $p$ is of maximal
Higgs field. Consequently one has the equality
$$c_1({E^{0,1}_{\phi}})=\frac{1}{2}c_1(\calM_0).$$  Analogous statements hold for the bar counterpart.
\end{proposition}
\begin{proof}
The Higgs subbundle $(E_{\phi}=E^{1,0}_{\phi}\oplus
E^{0,1}_{\phi},\theta_{\phi})$ is the modulo $p$ reduction of the
Higgs bundle $(\tilde E_{\phi}, \tilde \theta_{\phi})$ over $\calM$
by Corollary \ref{decomposition in char 0}. For $\phi|_{F}=\tau$,
the Higgs bundle $(E^0_{\phi}, \theta^0_{\phi})$, that is the
base change of $(\tilde E_{\phi}, \tilde \theta_{\phi})$ to
$\calM^0$, is actually of maximal Higgs field (see \cite{VZ}). That
is, the Higgs field ${\theta^0_{\phi}}^{1,0}: {E^0_{\phi}}^{1,0}\to
{E^0_{\phi}}^{0,1}\otimes \Omega_{\calM^0}$ is an isomorphism. Then
under the assumption on $p$, we claim that the Higgs field in char
$p$ must be maximal. In fact, the Higgs field
$\theta_{\phi}^{1,0}$ can not be zero. Otherwise, the Higgs
subbundle $(E_{\phi}^{1,0},0)$ of $(E_{\phi},\theta_{\phi})$ is of
non-positive degree by the Higgs semistability (see Proposition
\ref{app-Ogus-V}). This is in contradiction with the fact that
$$
\deg
E_{\phi}^{1,0}=\deg {E^0_{\phi}}^{1,0}=\frac{1}{2}\deg
\Omega_{\calM^0}=g-1>0.
$$
Since $\theta_{\phi}^{1,0}$ is a nonzero morphism between line
bundles with the same degree, it must be an isomorphism.\\

Moreover ${E^0_{\phi}}^{0,1}$ is isomorphic to the dual
${E^0_{\phi}}^{1,0*}$ of ${E^0_{\phi}}^{1,0}$. By the
theorem of Langton (Main Theorem A' in \cite{La}), the isomorphism
extends to an isomorphism $E_{\phi}^{0,1}\cong E_{\phi}^{1,0*}$.
Thus the maximality of $\theta_{\phi}$ implies that
$(E_{\phi}^{0,1})^{2}\cong \calT_{\calM_0}$. Taking the cycle classes
of both sides of the isomorphism, we obtain the claimed formula.
\end{proof}

We believe that one can remove the condition on the prime $p$ in the
above proposition. It is clear that these Higgs
subbundles in the decomposition of $(E,\theta)$ are divided into two
types: one is of maximal Higgs field and the other is of
zero Higgs field. This is the char $p$ analogue of the
corresponding result in \cite{VZ} in the char 0 case.

\section{The Newton Polygon Jumping Locus of the Shimura
Curve}\label{Newton jumping locus}

In this section we prove the mass formula for the Shimura curve
$\calM_0$. We refer to \cite{FL} for the definition of the Newton
polygon stratifications and other related notions. By Corollary
\ref{slope}, there are only two possible Newton polygons for points in
$\calM_0(\bbF)$. We denote by $\calS$ the subset of
$\calM_0(\bbF)$ consisting of the closed points for which the Newton polygon jumps.
By a theorem of Grothendieck-Katz (see
\cite{Ka1}), the Newton polygon jumps under specialization, and
$\calS$ is an algebraically closed subset of $\calM(\bbF)$.
In particular, the cardinality  of $\calS$ is finite. \\

We find that the
morphisms $\calF_{\calX_0|\calM_0}^{n}: F_{\calM_0}^{n*}E^{0,1}\to
E^{0,1}, \ n\geq 1$, where $\calF_{\calX_0|\calM_0}^{n}$ is the
composition of relative Frobenius morphisms (see \cite{Ka1} and
\cite{LSZ}), can be applied to compute the number $|\calS|$, as is very interesting.
%In our special situation, it is best to look into the
%restriction of $\calF_{\calX_0|\calM_0}^{n}$ to each direct summand
%in $E^{0,1}$ in Proposition \ref{decomposition}.
One notices that the restriction of $\calF_{\calX_0|\calM_0}$
induces a morphism $\calF_{\calX_0|\calM_0}:
F_{\calM_0}^{*}E^{0,1}_{\phi} \to E^{0,1}_{\sigma\phi}$ for each
$\phi\in \Phi$, since the Frobenius morphism on $\bbD$ is $\sigma$-semilinear.
Since each prime of $F$ is inert in $L$, we shall
use the same letter $\mathfrak{p}$ to denote the prime of $L$ lying over
the prime $\mathfrak{p}$ of $F$. We write the subset of $\Phi$ as
$$
Hom_{\bbQp}(L_{\frakp},\bbQpp)=\{\phi_1,\cdots,\phi_d,\phi_1^*,\cdots,\phi_d^*\},
$$
in such a way that $\phi_1|_F=\phi_1^*|_F=\tau$, and the
Frobenius automorphism $\sigma$, which is the generator of
$Gal(L_{\frakp}|\bbQ_p)$, acts on the set as the cyclic
permutation of $2d$ letters. For example, $\sigma \phi_d= \phi_1^*$,
$\sigma \phi_d^*= \phi_1$, and so on.
\begin{proposition}\label{isomorphisms between direct summands of de Rham bundle}
The notations are as above and all morphisms in the following are
the relative Frobenius morphisms. Let $\phi\in \Phi$. Then the
following statements hold:
\begin{itemize}
  \item [(i)] For $\phi\in Hom_{\bbQp}(L_{\frakp},\bbQpp)$, one has
  two possibilities:
  \begin{itemize}
    \item [(i.1)] If $d=1$, then $F_{\calM_0}^{*}E^{0,1}_{\phi_1}\to E^{0,1}_{\phi_1^*}$ is nonzero.
The same holds for $F_{\calM_0}^{*}E^{0,1}_{\phi_1^*}\to
E^{0,1}_{\phi_1}$ and the bar counterparts. Moreover, they have the
same zero locus.
    \item [(i.2)] If $d\geq 2$, then $F_{\calM_0}^{*}E^{0,1}_{\phi_1}\to
E^{0,1}_{\phi_2}$ is injective, $F_{\calM_0}^{*}E^{0,1}_{\phi_d}\to
E^{0,1}_{\phi_1^*}$ is surjective, and for the bar counterparts,
$$(\calF_{\calX_0|\calM_0})|_{F_{\calM_0}^{*}E^{0,1}_{\bar{\phi}_1}}=(\calF_{\calX_0|\calM_0})|_{F_{\calM_0}^{*}E^{0,1}_{\bar{\phi}_1^*}}=0.$$
Moreover if $d\geq 3$, then for $2\leq i\leq d-1$,
$F_{\calM_0}^{*}E^{0,1}_{\phi_i}\to E^{0,1}_{\phi_{i+1}}$ is an
isomorphism.
\end{itemize}
  \item [(ii)] For $\phi\notin Hom_{\bbQp}(L_{\frakp},\bbQpp)$, one has an isomorphism
  $$
\calF_{\calX_0|\calM_0}: F_{\calM_0}^{*}E^{0,1}_{\phi}\cong
E^{0,1}_{\sigma\phi}.
$$
\end{itemize}
\end{proposition}
\begin{proof}
We adopt a pointwise argument here. Let $t\in \calM_0(\bbF)$ and let
$A$ be the fiber of $f_0$ over $t$. Let $(\bbD(A),\calF,\calV)$ be
the Dieudonn\'{e} module of $A$. We first discuss the $d\geq 2$
case. Since
$\calF_{\calX_0|\calM_0}({F_{\calM_0}^{*}E^{0,1}_{\bar{\phi}_1}})\subset
E^{0,1}_{\sigma\bar{\phi}_1}$ and $E^{0,1}_{\sigma\bar{\phi}_1}=0$
for $(\sigma\bar{\phi}_1)|_F\neq \bar\tau$, we have
$\calF_{\calX_0|\calM_0}|_{F_{\calM_0}^{*}E^{0,1}_{\bar{\phi}_1}}=0$
and similarly for the $\bar{\phi}_1^*$-summand. Now we assume
$d\geq 3$ and look at the morphism
$F_{\calM_0}^{*}E^{0,1}_{\phi_i}\cong E^{0,1}_{\phi_{i+1}}$ for
$2\leq i\leq d-1$. Without loss of generality we discuss this only for $i=2$.
Since $\phi_i|_F\neq \tau$ for $i=2,3$, we have
$\dim_{\bbF}\frac{\bbD(A)_{\phi_i}}{p\bbD(A)_{\phi_i}}=2$.
Furthermore $\calV(\bbD(A)_{\phi_3})=p\bbD(A)_{\phi_2}$ by
Proposition \ref{Dieudonne module structure} (iii). We consider the
$\sigma$-semilinear map
$$
\calF \mod p: \frac{\bbD(A)_{\phi_2}}{\calV\bbD(A)_{\phi_2}}\to
\frac{\bbD(A)_{\phi_3}}{\calV\bbD(A)_{\phi_3}}
$$
induced by $\calF: \bbD(A)_{\phi_2}\to \bbD(A)_{\phi_3}$. It is
known that the above map is simply the Hasse-Witt map after
identifying the spaces
$\frac{\bbD(A)_{\phi_i}}{p\bbD(A)_{\phi_i}}=H^1(A,\calO_A)_{\phi_i}$.
Now for $e\in \bbD(A)_{\phi_2}$, $e\mod p\in \ker(\calF \mod p)$ if and only if
$\calF(e)\in p\bbD(A)_{\phi_3}$, and if and only if $e\in
\calV(\bbD(A)_{\phi_3})=p\bbD(A)_{\phi_2}$ (by applying $\calV$ or
$\calF$ to both sides). So one sees that $\calF \mod p$ is injective.
Hence it must be surjective for the dimensional reason. This
proves the isomorphism in (i.2) for $i=2$ and similarly we have the
isomorphisms in (ii). The similar argument proves the
injectivity of $F_{\calM_0}^{*}E^{0,1}_{\phi_1}\to E^{0,1}_{\phi_2}$
and the surjectivity of
$F_{\calM_0}^{*}E^{0,1}_{\phi_d}\to E^{0,1}_{\phi_1^*}$ follows by duality. It remains
to show the $d=1$ case. In this case, the result in Corollary
\ref{slope} tells us the $p$-rank of abelian varieties in
$\calM_0(\bbF)$ is either $4(n-1)$ or $4n$. Now by (ii), all of the
$\phi$-summands with $\phi\notin Hom_{\bbQp}(L_{\frakp},\bbQpp)$
have contributions to the $p$-rank. The existence of
closed points of $p$-rank $4(n-1)$ implies that each of the four
morphisms in (i) can not be zero. Moreover, because of the
existence of the other $p$-rank, all of the four
morphisms will vanish at a point $t$ as soon as one of them vanishes at $t$.
\end{proof}

Now we consider the composition of the relative Frobenius morphisms
$$
{\calF_{\calX_0|\calM_0}^{d}}: {E^{0,1}_{\phi_1}}^{(p^d)}\rightarrow
{E^{0,1}_{\phi_2}}^{(p^{d-1})} \rightarrow\cdots\rightarrow
{E^{0,1}_{\phi_d}}^{(p)}\rightarrow E^{0,1}_{\phi_1^*}.
$$
As a consequence of the above analysis of the relative Frobenius
morphisms, we have the following
\begin{proposition}\label{Hasse-Witt locus}
Let $t\in \calM_0(\bbF)$. Then $t\in \calS$ if and only if the map
${\calF_{\calX_0|\calM_0}^{d}}: {E^{0,1}_{\phi_1}}^{(p^d)}\to
E^{0,1}_{\phi_1^*}$ vanishes at $t$.
\end{proposition}

\begin{proof}
For simplicity we put $h=\calF_{\calX_0|\calM_0}^{d}$. The $d=1$
case follows directly from Proposition \ref{isomorphisms between
direct summands of de Rham bundle} (i). Assume $d\geq 2$. One
direction is clear. Namely if $t\in S$, then the restriction of $h$
to the $\phi_1$-summand must be zero at $t$ by Proposition
\ref{isomorphisms between direct summands of de Rham bundle} (ii). We
proceed to show the converse direction. Again by Proposition
\ref{isomorphisms between direct summands of de Rham bundle} (ii) the
$\phi$-summands with $\phi\notin Hom_{\bbQp}(L_{\frakp},\bbQpp)$ have no contribution
to the jumping of the Newton polygons and we only need
to consider $\phi\in Hom_{\bbQp}(L_{\frakp},\bbQpp)$. It is clear
that $h$ maps the $\phi_i$-summand to the $\phi_i^*$-summand and the
$\phi_i^*$-summand to the $\phi_i$-summand for $1\leq i\leq d$. For
example, one has $h({E^{0,1}_{\phi_1^*}}^{(p^d)})\subset
E^{0,1}_{\phi_1}$. Note that for the $\phi$-summand with $\phi|_F\neq
\tau$ the rank of $h$ must be reduced by at least one and thus it contributes
to the $p$-rank at most by one. There are $2(d-1)$ such summands in
all. By the assumption $h$ vanishes at $t$ on the $\phi_1$-summand.
Applying the endomorphism $\Pi$ one sees that $h$ vanishes at $t$ on
the $\phi_1^*$-summand too. It implies that the $p$-rank at $t$ is at most
$4n-2d-2$ by Proposition \ref{isomorphisms between direct
summands of de Rham bundle} (i.2), and therefore the $p$-rank can not be $4n-2d$.
Since there are only two possibilities for the $p$-ranks at $t\in S$.
From the above proof, we also see that the restriction of $h$ to each
$\phi$-summand with $\phi\in Hom_{\bbQp}(L_{\frakp},\bbQpp)$ vanishes
at $t$ if and only if $t\in S$.
\end{proof}

Furthermore, we have
\begin{proposition}\label{multiplicity free}
The zero locus of ${\calF_{\calX_0|\calM_0}^{d}}:
{E^{0,1}_{\phi_1}}^{(p^d)}\to E^{0,1}_{\phi_1^*}$ is reduced. In
other words, the zero divisor of the global section of the line
bundle $({E^{0,1}_{\phi_1}}^{(p^d)})^{-1}\otimes E^{0,1}_{\phi_1^*}$
defined by ${\calF_{\calX_0|\calM_0}^{d}}$ is of multiplicity one.
\end{proposition}
Before we show this result, we make a digression into the display theory of
Dieudonn\'{e} modules. By a theorem of Serre-Tate
(\cite{Ka2}), the equi-characteristic deformation of an abelian
variety $A$ in positive characteristic is the same as that of its
$p$-divisible group $A(p)$. The latter is determined by the display
of the Dieudonn\'{e} module $\bbD(A^{t}(p))$ (see \cite{NO} and the
references therein). This is also true when polarizations and
endomorphisms are considered. We put $\bbD=\bbD(A^{t}(p))$ and recall
that we have the following decomposition:
\begin{eqnarray*}\bbD=\bigoplus_{\phi\in
\Phi}(\bbD_{\phi}\oplus \bbD_{\bar{\phi}}),\end{eqnarray*} where
$\Phi=Hom_{\bbQ}(L,\bar \bbQ)$.

\begin{lemma} \label{basis}
There exists a basis $\{X_{\phi},Y_{\phi},
X_{\bar{\phi}},Y_{\bar{\phi}}\ |\ \phi\in \Phi\}$ of $\bbD$, such
that\par (i)  $\{X_{\phi},Y_{\phi}\}$ is a basis for $\bbD_{\phi}$,
and $X_{\bar{\phi}}, Y_{\bar{\phi}}$ is the dual basis of
$\bbD_{\bar{\phi}}$.
\par (ii) For $\phi|_F=\tau$, $Y_{\phi}, X_{\bar{\phi}}\in \calV(\bbD)$;
while for $\phi|_F\neq\tau$, $X_{\bar{\phi}}, Y_{\bar{\phi}}\in
\calV(\bbD)$.
\end{lemma}
\begin{proof}
Note that for $\phi\in \Phi$, we have a short exact sequence:
$$0\rightarrow \calV\bbD_{\sigma\phi}/p\bbD_{\phi}\rightarrow\bbD_{\phi}/p\bbD_{\phi}
\rightarrow\bbD_{\phi}/\calV\bbD_{\sigma\phi}\rightarrow 0.$$ For
$\phi|_F\neq\tau$, the statements in (i) and (ii) are obvious, by
Proposition~\ref{Dieudonne module structure}. For $\phi|_F=\tau$,
first of all, we prove that we can choose a basis
$X_{\phi},Y_{\phi}$ for $\bbD_{\phi}$, such that $Y_{\phi}\in
\calV(\bbD_{\sigma\phi})$. In this case, by
Proposition~\ref{Dieudonne module structure}, the dimensions of the
terms appearing in the above exact sequence are in turn $1,2,1$. Let
$X_{\phi},Y_{\phi}\in \bbD_{\phi}$, such that the image of
$\bar{X}_{\phi}$ generates $\bbD_{\phi}/\calV\bbD_{\sigma\phi}$ as
vector space, and $\bar{Y}_{\phi}$ generates
$\calV(\bbD_{\sigma\phi})$ as vector space. Then
$\bar{X}_{\phi},\bar{Y}_{\phi}$ generate $\bbD_{\phi}/p\bbD_{\phi}$.
By Nakayama's Lemma, $X_{\phi},Y_{\phi}$ is a basis of $\bbD_{\phi}$
with $Y_{\phi}\in \calV(\bbD_{\sigma\phi})$. Secondly, we prove that
there is a dual basis $X_{\bar{\phi}},Y_{\bar{\phi}}$ of
$\bbD_{\bar{\phi}}$, such that $X_{\bar{\phi}}\in \calV(\bbD)$.
Similarly as above, we can find a basis $x,y$ of $\bbD_{\bar{\phi}}$
with $y\in \calV(\bbD)$. Let $H=\begin{pmatrix}
a   & b \\
c   & d
\end{pmatrix}$ be the intersection matrix of $X_{\phi},Y_{\phi}$ and $x,y$.
Thus we have the valuation $v_p(d)>0$, $b,c$ are invertible and $H$
is invertible. By solving a system of linear equations, we see that
$X_{\bar{\phi}}=\frac{1}{det(H)}(dx-cy)$ and
$Y_{\bar{\phi}}=\frac{1}{det(H)}(-bx+ay)$ satisfy the requirements,
since $v_p(d)>0$ implies that
$dx=pd'x=\calV\calF(d'x)\in\calV(\bbD)$.

\end{proof}

Under this basis of $\bbD$, the corresponding display is
$$
\begin{pmatrix}
A   & B \\
C   & D
\end{pmatrix},$$
where the matrix $A,C$ are (we take $n=2, d=[F_{\frakp}:\bbQp]=2$
for example)
$$A=\begin{pmatrix}
0&0&0&0&a_1&c_1&0&0\\
a_2&0&0&0&0&0&0&0\\
b_2&0&0&0&0&0&0&0\\
0&a_{1^*}&c_{1^*}&0&0&0&0&0\\
0&0&0&a_{2^*}&0&0&0&0\\
0&0&0&b_{2^*}&0&0&0&0\\
0&0&0&0&0&0&0&0\\
0&0&0&0&0&0&0&0
\end{pmatrix}$$
and $$C=\begin{pmatrix}
0&0&0&0&b_1&d_1&0&0\\
0&0&0&0&0&0&e_1&0\\
0&0&0&0&0&0&f_1&0\\
0&b_{1^*}&d_{1^*}&0&0&0&0&0\\
0&0&0&0&0&0&0&e_{1^*}\\
0&0&0&0&0&0&0&f_{1^*}\\
0&0&0&0&0&0&0&0\\
0&0&0&0&0&0&0&0
\end{pmatrix}.
$$
Here the basis is arranged in an obvious manner. \iffalse Namely,
$$
\{X_{\phi_1},X_{\phi_2},Y_{\phi_2},X_{\phi_{1^*}},X_{\phi_{2^*}},Y_{\phi_{2^*}},Y_{\bar{\phi}_{1}},Y_{\bar{\phi}_{1^*}}\}
$$
and
$$
\{Y_{\phi_1},X_{\bar{\phi}_{2}},Y_{\bar{\phi}_{2}},Y_{\phi_{1^*}},X_{\bar{\phi}_{2}^*},Y_{\bar{\phi}_{1}^*},X_{\bar{\phi}_{1}},X_{\bar{\phi}_{2}}\}.
$$
\fi
In this case, the Frobenius is given by the matrix $\begin{pmatrix}A&pB\\
C&pD\end{pmatrix}$.

\begin{lemma}
Let $R=k[[t]]$. The $\calO_B$-action can be extended to $spec (R)$.
Hence the display of the infinitesimal universal deformation is
given by
$$\begin{pmatrix}A+TC&p(B+TD)\\ C&pD\end{pmatrix},$$
where $T$ is the Teichm\"{u}ller lifting of $t$ (that is,
$T=(t,0,\cdots)$). In particular, the matrix $A+TC$, read mod $p$,
is the Hasse-Witt matrix of the deformation corresponding to $T$.
\end{lemma}
\begin{proof}
It is known from Proposition~\ref{local-ring} that the local
deformation ring of the Shimura curve is regular on one parameter.
Let $\bbD_R$ be the display over $R$. Then
$$\bbD_R=\bigoplus_{\phi\in \Phi} (\bbD_{R,\phi}\bigoplus \bbD_{R,\bar{\phi}}),$$
where $\bbD_{R,\phi}$ (resp. $\bbD_{R,\bar{\phi}}$) is obtained from
$\bbD_{\phi}$ (resp. $\bbD_{\bar{\phi}}$) by extending scalars to
$W(R)$, with the naturally given action of $W(k)$ on each component.
Recall that the action of $\calO_{LK}$ is defined via the map
$$\calO_{LK}\rightarrow \bigoplus_{\phi}
(W(k)\oplus W(k)), a\mapsto
(\cdots,w(\phi)(a),\bar{w}(\phi)(a),\cdots).$$ This is a map of
Dieudonn\'{e} modules if and only if it commutes with the Frobenius; that is,
if and only if
$$M_1M_2^{\sigma}=M_2M_1,$$
where $M_1=\begin{pmatrix}A+TC&p(B+TD)\\
C&pD\end{pmatrix}$ and $$M_2=\begin{pmatrix}diag(\cdots,w(\phi)(a),\cdots,\bar{w}(\phi)(a),\cdots)& 0\\
0& diag(\cdots,w(\phi)(a),\cdots,\bar{w}(\phi)(a),\cdots)
\end{pmatrix}.$$ It is easy to verify that this is true, by a direct computation.
\end{proof}

\iffalse According to Corollary~\ref{slope}, we know that there are
only two different Newton polygons $P_1,P_0$ for points in
$\calM(k)$, with $P_0$ lies above $P_1$ (we denote this by $P_0\prec
P_1$).
\begin{theorem}\label{main1}
Let $t_0\in \calM(\bbF)$ be a geometric point with Newton polygon
$P_{0}$. There is a choice of parameter $t$, such that the universal
local deformation ring of $X_{t_0}$ is isomorphic to $k[[t]]$, and
the closed regular formal subscheme which parameterizes the abelian
varieties $\calX_s$ with $P_s= P_{0}$ is defined by the ideal $(t)$,
and the multiplicity is one.
\end{theorem}
\fi

Now we come to the proof of Proposition \ref{multiplicity free}.
\begin{proof}
The universal Dieudonn\'{e} module $\bbD_R$ is displayed by the
matrix
$$\begin{pmatrix}
A+TC& B+TD\\
C&D
\end{pmatrix}.$$
Now we show that the locus of $\calF^d$ is of multiplicity one. For
this we put
$$
\calV^d=\calV_{\calX_0|\calM_0}^d: E^{1,0}_{\bar \phi_1^*}\to
{E^{0,1}_{\bar \phi_1}}^{p^d},
$$
which is the dual of $\calF^d$. It is then equivalent to show the zero locus
of $\calV^d$ is of multiplicity one. Moreover, for
simplicity we just take $n=d=2$ in the following argument, and the argument
for the general case is completely the same. The following matrix
mod $p$ is the Hasse-Witt matrix of the Frobenius of the local
deformation:
$$A+TC=\begin{pmatrix}
0&0&0&0&a_1+tb_1&c_1+td_1&0&0\\
a_2&0&0&0&0&0&te_1&0\\
b_2&0&0&0&0&0&tf_1&0\\
0&a_{1^*}+tb_{1^*}&c_{1^*}+td_{1^*}&0&0&0&0&0\\
0&0&0&a_{2^*}&0&0&0&te_{1^*}\\
0&0&0&b_{2^*}&0&0&0&tf_{1^*}\\
0&0&0&0&0&0&0&0\\
0&0&0&0&0&0&0&0
\end{pmatrix}.$$
Thus the Hasse-Witt matrix of $\calV^2$ is the following matrix mod
$p$:
$$(A+TC)(A+TC)^{\sigma}=\begin{pmatrix}
0&0&0&f_{14}&0& 0&0&0\\
0&0&0&0&f_{25}&f_{26}&0&0\\
0&0&0&0&f_{35}&f_{36}&0&0\\
f_{41}&0&0&0&0&0&0&0\\
0&f_{52}&f_{53}&0&0&0&0&0\\
0&f_{62}&f_{63}&0&0&0&0&0\\
0&0&0&0&0&0&0&0\\
0&0&0&0&0&0&0&0
\end{pmatrix},$$
where
\begin{eqnarray*}
% \nonumber to remove numbering (before each equation)
f_{14}&=&(a_1+tb_1)a_{2^*}^{\sigma}+(c_1+td_1)b_{2^*}^{\sigma}, \\
f_{25}&=&a_2(a_1^{\sigma}+t^{\sigma}b_1^{\sigma}),f_{26}=a_2(c_1^{\sigma}+t^{\sigma}d_1^{\sigma}),\\
f_{35}&=&b_2(a_1^{\sigma}+t^{\sigma}b_1^{\sigma}), f_{36}=b_2(a_1^{\sigma}+t^{\sigma}b_1^{\sigma}),\\
f_{41}&=&(a_{1^*}+tb_{1^*})a_{2}^{\sigma}+(c_{1^*}+td_{1^*})b_{2}^{\sigma},\\
f_{52}&=&a_{2^*}(a_{1^*}^{\sigma}+t^{\sigma}b_{1^*}^{\sigma}),f_{53}=a_{2^*}(c_{1^*}^{\sigma}+t^{\sigma}d_{1^*}^{\sigma}),\\
f_{62}&=&b_{2^*}(a_{1^*}^{\sigma}+t^{\sigma}b_{1^*}^{\sigma}),\\
f_{63}&=&b_{2^*}(c_{1^*}^{\sigma}+t^{\sigma}d_{1^*}^{\sigma}).
\end{eqnarray*}
Thus $\calV^2$ is locally given by the function
$$(a_{1^*}+tb_{1^*})a_{2}^{\sigma}+(c_{1^*}+td_{1^*})b_{2}^{\sigma}\\
=(a_{1^*}a_{2}^{\sigma}+c_{1^*}b_{2}^{\sigma})+t(b_{1^*}a_{2}^{\sigma}+d_{1^*}b_{2}^{\sigma})=0.$$
As $t=0\in \calS$, we have
$a_{1^*}a_{2}^{\sigma}+c_{1^*}b_{2}^{\sigma}=0$ and
$b_{1^*}a_{2}^{\sigma}+d_{1^*}b_{2}^{\sigma}\neq 0$. Thus the locus
is of multiplicity one. \hspace*{9.5cm}
\end{proof}

\begin{theorem}\label{mass formula}
Let $\calM_0$ be the moduli space constructed in \S2 and let $\calS$
be the Newton polygon jumping locus in $\calM_0$. Then in the Chow
ring of $\calM_0$ the following formula holds:
$$
\calS= \frac{1}{2}(1-p^d)c_1(\calM_0).
$$
where $d=[F_{\mathfrak{p}}:\bbQ_p]$ is the local degree. As a
consequence, one obtains the mass formula for $\calM_0$:
$$
|\calS|=(p^d-1)(g-1),
$$
where $g$ is the genus of the Shimura curve $\calM_0$.
\end{theorem}
\begin{proof}
The mass formula follows by taking the degree in the cycle formula.
By Proposition \ref{uniformizing Higgs bundle} and \ref{maximal Higgs field in char p}, the
cycle of the zero locus of ${F_{\calX_0|\calM_0}^{d}}:
{E^{0,1}_{\phi_1}}^{(p^d)}\to E^{0,1}_{\phi_1^*}$ is equal to
$$
c_1(E^{0,1}_{\phi_1^*})-c_1({E^{0,1}_{\phi_1}}^{(p^d)})=
(1-p^d)c_1({E^{0,1}_{\phi_1}})= \frac{1}{2}(1-p^d)c_1(\calM_0).
$$
Then the theorem follows from Proposition \ref{Hasse-Witt locus} and
\ref{multiplicity free}.
\end{proof}
\iffalse
\begin{remark}
Even if $C^p$ in the definition of the level structure $C=C_p\times
C^p\subset G(\bbA_f)$ (cf. \S2) is so large that the resulting
coarse moduli space $\calM_0$ supports no universal family, the mass
formula still holds for $\calM_0$. One can always take a finite
Galois \'{e}tale covering $\tilde \calM_0$ of $\calM_0$ with large
degree so that it supports a universal family. Since both sides of
the formula increases by multiplying the degree of the covering
under \'{e}tale base change, one obtains the formula for $\calM_0$
from that for $\tilde \calM_0$.
\end{remark}
\fi
\begin{remark}
The last two sections have certain overlaps with parts of the
paper \cite{Kas} by P. Kassaei. In particular one shall compare
Corollary \ref{maximal Higgs field in char p} and Proposition
\ref{multiplicity free} with Proposition 4.1 and 4.3 in \cite{Kas}.
\end{remark}

\section{Stability and instability of the Higgs subbundles in char $p$}\label{stability}
In this section we study the stability and instability of
the Higgs subbundles constructed in Proposition \ref{decomposition}.
We will assume that $p\geq 2g$ in this section, unless otherwise specified.
\begin{proposition}\label{app-Ogus-V}
 With the assumption on $p$ as above, we have that for each $\phi\in\Phi$,
$(E_{\phi},\theta_{\phi})$ and $(E_{\bar \phi},\theta_{\bar \phi})$
are Higgs-semistable of degree $0$. Particularly for
$\phi|_F\neq\tau$, the rank two vector bundles $E^{0,1}_{\phi}$ and
$E^{1,0}_{\bar\phi}$ are semistable.
\end{proposition}
\begin{proof}
By construction, for each $\phi\in \Phi$, $(E_{\phi},\theta_{\phi})$
and $(E_{\bar \phi},\theta_{\bar \phi})$ are the modulo $p$
reductions of Higgs bundles in characteristic 0 by Corollary
\ref{decomposition in char 0}. By Theorem 4.14 (3) and Proposition
4.19 in \cite{OV}, they are Higgs-semistable under the assumption
on $p$ as in the statement. Moreover for each place $\phi$ with
$\phi|_F\neq\tau$ one has
$(E_{\phi},\theta_{\phi})=(E_{\phi}^{0,1},\theta_{\phi}^{0,1})$ and
$(E_{\bar \phi},\theta_{\bar \phi})=(E_{\bar
\phi}^{1,0},\theta_{\bar \phi}^{1,0})$ by Proposition
\ref{decomposition}. The Higgs field $\theta_{\phi}^{0,1}$ is by
definition zero, and $\theta_{\bar \phi}^{1,0}$ is also zero as
$E_{\bar \phi}^{0,1}$ is a zero bundle.
\end{proof}
%In the following discussion the assumption on $p$ in the above
%proposition will be assumed.
For a semistable bundle $E$ of rank
$r$ over a smooth projective curve $C$ defined over $\bbF$, one can
ask further the semi-stability of the bundle over $C$ under the $n$-th
iterated Frobenius pull-back $F_{C}^{*n}E$ for $n\geq 1$. It turns
out that the bundles $F_{C}^{*n}E$ are not necessarily semistable.
In order to measure the instability of $F_C^*E$ one introduces and
studies the invariant
$\nu(F_C^*E)=\mu_{max}(F_{C}^*E)-\mu_{min}(F_{C}^*E)$ where
$\mu_{max}(F_{C}^*E)$ (resp. $\mu_{min}(F_{C}^*E)$) is the slope of
$E_1$ (resp. $\frac{E_n}{E_{n-1}}$) in the Harder-Narasimhan
filtration of $F_C^*E$:
$$
0=E_0\subset E_1\subset \cdots\subset E_n=F_C^*E.
$$
The following result says that $F_C^*E$ can not be very instable by
exhibiting an upper bound of $\nu(F_C^*E)$.
\begin{theorem}[Lange-Stuhler, Satz 2.4 \cite{LS} for $r=2$; Shepherd-Barron, Corollary 2 \cite{SB} and Sun, Theorem 3.1 \cite{Sun1} for arbitrary rank]\label{instability}
Let $E$ be a rank $r$ semistable bundle over a smooth projective
curve $C$ of genus $g$ defined over $\bbF$. Then one has the
inequality
$$
\nu(F_C^*E)\leq (r-1)(2g-2).
$$
In particular, $F_C^*E$ is still semistable when $g\leq 1$.
\end{theorem}
Based on the above inequality one can also deduce a generalization
of it for $n\geq 2$ (see Theorem 3.7 in \cite{LSZ}). A. Langer
(Corollary 6.2 in \cite{LG}) actually obtained a better bound on
this issue and generalized the above inequality as well to a higher
dimensional base. It is then interesting to find examples where
the upper bound of the inequality is reached. In the case where the local
degree $[F_{\frakp}:\bbQp]$ is strictly larger than one, certain
Higgs subbundles over the Shimura curve $\calM_0$ do provide such
examples (See the proposition of \S4.4 in  \cite{JRXY} for a
classification of semistable bundles of rank two over curves in char
$2$).
\begin{proposition}\label{instable}
Assume that $[F_{\frakp}:\bbQp]>1$. The rank two semistable bundles
$E^{0,1}_{\phi_d}$ and $E^{0,1}_{\phi_d^*}$ over $\calM_0$ achieve
the upper bound in Theorem \ref{instability}. That is, one has the
equality
$\nu({E^{0,1}_{\phi_d}}^{(p)})=\nu({E^{0,1}_{\phi_d^*}}^{(p)})=2g-2$.
\end{proposition}
\begin{proof} It suffices to show the result for the $\phi_d$ piece.
The proof for another piece is completely similar.
We consider the morphism
$$
F_{\calX_0|\calM_0}: {E^{0,1}_{\phi_d}}^{(p)}\rightarrow
E^{0,1}_{\phi_1^*},
$$
where $E^{0,1}_{\phi_1^*}$ is a line bundle since $\phi_1^*|_F=\tau$. We put then $\calE$ (resp. ${\calE}'$)
to be the kernel (resp. the image) of the above morphism. That is,
we have the following short exact sequence:
$$
0\rightarrow \calE \rightarrow  {E^{0,1}_{\phi_d}}^{(p)}\rightarrow
\calE'\rightarrow 0.
$$
Since $\deg({E^{0,1}_{\phi_d}})=0$ and
$\deg(E^{0,1}_{\phi_1^*})=1-g$, we have $\deg {E^{0,1}_{\phi_d}}^{(p)}=0$
and $\deg \calE'\leq 1-g$. So by the inequality in Theorem
\ref{instability}, one has the following inequalities
$$
\mu(\calE)-\mu(\calE')\leq \nu({E^{0,1}_{\phi_d}}^{(p)})\leq
2g-2\leq \mu(\calE)-\mu(\calE').
$$
It follows that the above inequalities have to be an equality at each
step and particularly the assertion of the proposition follows.
\end{proof}
Combining this proposition with Proposition \ref{isomorphisms between direct summands of
de Rham bundle} (ii), one obtains Theorem \ref{thm 3} (iii). \\

The other extreme about the semi-stability of a vector bundle under
iterated Frobenius pull-backs is expressed in the following
definition.
\begin{definition}
Let $E$ be a semistable vector bundle over a smooth projective
curve $C$ as above. It is said to be strongly semistable if for all
$n\geq 1$, the $F_C^{*n}E$ are semistable.
\end{definition}
In the case that for certain $n\geq 1$ one has an isomorphism
$F_C^{*n}E\cong E$, the bundle $E$ is obviously strongly
semistable. The following theorem gives us a characterization of the
category of the strongly semistable bundles over $C$.
\begin{theorem}[Lange-Stuhler, \S1 \cite{LS}]\label{result of LS}
Let $C$ be a smooth projective curve as above and $\pi_1^{et}(C)$ be
its \'{e}tale fundamental group. Let $E$ be a vector bundle over $C$
of rank $r$. Then the following conditions about $E$ are equivalent:
\begin{itemize}
  \item [(i)]  There exists $n\geq 1$ such that $F_C^{*n}E\cong E$,
  \item [(ii)] There exists an \'{e}tale covering map $\pi: \tilde C\to C$ such that $\pi^*E$ over $\tilde C$ is trivial,
  \item [(iii)] $E$ corresponds to a continuous representation $\pi_1^{et}(C)\to Gl_{r}(\bbF)$ where $Gl_{r}(\bbF)$ is equipped with the discrete topology.
\end{itemize}
One calls the bundle $E$ satisfying one of the above equivalent
conditions \'{e}tale trivializable. The bundle $E$ is strongly
semistable if and only if there exists $n\geq 0$ such that $F_C^{*n}E$ is
\'{e}tale trivializable.
\end{theorem}
We can find such examples again among the Higgs subbundles in this study.
\begin{proposition}\label{sss}
%\begin{itemize}
  (i) Assume $p$ is not inert in $F$. Then for $\phi\notin Hom_{\bbQp}(L_{\frakp},\bbQpp)$,
$E^{0,1}_{\phi}$ and $E^{1,0}_{\bar \phi}$ are \'{e}tale
trivializable, and particularly strongly semistable.   \par

 (ii) In case $d=[F_{\frakp}:\bbQp]>1$, the semistable bundles $E^{0,1}_{\phi_i}$ for $2\leq i\leq d$ are
not strongly semistable and consequently stable.
%\end{itemize}
\end{proposition}
\begin{proof}
(i)  By Proposition \ref{isomorphisms between direct summands of
de Rham bundle} (ii), the morphism
$\calF_{\calX_0|\calM_0}:{E^{0,1}_{\phi}}^{(p)}\to
E^{0,1}_{\sigma\phi}$ is an isomorphism in the case that $\phi|_F\neq
\tau$. When $\phi\notin Hom_{\bbQp}(L_{\frakp_1},\bbQpp)$,
$(\sigma \phi)|_F\neq \tau$. This implies that for $n$ large enough
the composition of the relative Frobenious morphisms induces an
isomorphism ${E^{0,1}_{\phi}}^{(p^n)}\cong E^{0,1}_{\phi}$. The
proof for $E^{1,0}_{\bar
\phi}$ is similar by replacing the relative Frobenius morphism in the argument
by the relative Verschiebung morphism $V_{\calX_0|\calM_0}$. \\

(ii) For $\phi_i\in Hom_{\bbQp}(L_{\frakp_1},\bbQpp)$ with $i\neq
1$, we consider the composition $\eta$ of the morphisms
$$
{E^{0,1}_{\phi_i}}^{(p^{d-i+1})}\stackrel{F_{\calX_0|\calM_0}}{\cong}{E^{0,1}_{\phi_{i+1}}}^{(p^{d-i})}\cdots\stackrel{F_{\calX_0|\calM_0}}{\cong}
{E^{0,1}_{\phi_d}}^{(p)} \twoheadrightarrow E^{0,1}_{\phi_1^*},
$$
where the last morphism is surjective by the proof of
Proposition \ref{instable}. Thus the kernel of $\eta$ provides a sub line bundle of positive degree
(actually it is equal to $g-1$) in ${E^{0,1}_{\phi_i}}^{(p^{d-i+1})}$ and
therefore it is not semistable. The assertion about the stability
follows from the following simple lemma.
\end{proof}

\blem A strictly semistable rank two vector bundle $E$ of degree
zero over a smooth projective curve $C$ over $\bbF$ is strongly
semistable \elem
\begin{proof}
By the assumption, $E$ is an extension of two degree zero line
bundles. We write
$$
0\to \calL_1\to E\to \calL_2\to 0.
$$
Suppose that $E^{(p)}$ is not semistable, and let $\calL$ be a
positive sub line bundle of it. Then we have the morphisms
$$
\calL\to E^{(p)} \to \calL_2^{(p)}.
$$
Because $\deg \calL^{(p)}_2=0<\deg \calL$, the above composition is
zero. Hence the inclusion of $\calL$ in $E^{(p)}$ factors through
$\calL_1^{(p)}$. That is, one has $\calL \subset \calL_1^{(p)}$.
Again because $\deg \calL_1^{(p)}=0<\deg \calL$, we obtain a
contradiction. In conclusion, $E^{(p)}$ is semistable. By induction
on the iterations of Frobenius pull-backs we see that $E$ is actually
strongly semistable.
\end{proof}

\section{The Simpson-Ogus-Vologodsky correspondence of the Higgs subbundles in char
$p$}\label{SOV correspondence for Higgs subbundles} One of the main
results in \cite{OV} is to establish a char $p$ analogue of the
Simpson correspondence. Let $\pi:\calM_0'\to \calM_0$ be the
projection map and $\calF_{\calM_0}: \calM_0\to \calM_0'$ be the
relative Frobenius in the commutative diagram of Frobenius morphisms
for $\calM_0$ over $\mathbb F$ (see Notations and Conventions in \S1).
For the Shimura curve $\calM_0$ over $\mathbb F$, which is the
reduction of a Shimura curve over mixed characteristic, the Cartier
transform $C_{\calM_0}$ (see \cite{OV}) is a functor from the
category $MIC(\calM_0)$ of flat bundles over $\calM_0$ to the
category of Higgs bundles $HIG(\calM_0')$ over $\calM_0'$, with each of them
subject to suitable nilpotence conditions. The functor is an
equivalence of categories with the quasi-inverse functor $C^{-1}_{\calM_0}: HIG(\calM_0')\to
MIC(\calM_0)$. For the full subcategory of
flat bundles with vanishing $p$-curvatures, the functor is just the
classical Cartier descent (see \S5 in \cite{Ka0}), and it maps onto
the full subcategory of Higgs bundles with trivial Higgs fields over
$\calM_0'$. It is also clear that the functor transforms the
relative de Rham bundle $(\calH^1_{dR},\nabla)$ of the universal
family $f_0: \calX_0\to \calM_0$ to
$(E',\theta')=\pi_{}^{*}(E,\theta)$, which is the associated Higgs
bundle of the family $f_0':\calX_0'\to \calM_0'$, where $f_0'$ is the base change of
$f_0$ via $\pi$. In this section we examine the
Simpson-Ogus-Vologodsky correspondence and the Cartier transform for
the Higgs subbundles in Proposition \ref{decomposition}. It will be
assumed in this section that $p\geq \max\{2g,2(n+1)\}$, in order to
fulfill the basic requirements in Theorem 3.8, \cite{OV}. We use
$\calP'$ to denote the pull-back to $\calM_0'$
via $\pi$ of an algebra-geometric object $\calP$ defined over $\calM_0$. \\

We have the ring homomorphisms
$\calO_{LK}\stackrel{\varsigma}{\longrightarrow}
End_{\calM_0}(\calX_0)\stackrel{'}{\longrightarrow}End_{\calM_0'}(\calX_0')$.
Thus $\varsigma'$ denotes the composition morphism.
By abuse of notations, the induced ring homomorphisms from
$\calO_{LK}$ to the endomorphism rings $End(\calH^1_{dR},\nabla)$ and
$End(E,\theta)$ respectively are again denoted by $\varsigma$. For example, the
$\calO_{LK}$-action on $(E',\theta')$ via $\varsigma'$ decomposes
it into the direct sum of eigen Higgs subbundles
$$
(E',\theta')=\bigoplus_{\phi\in
\Phi}(E'_{\phi},\theta'_{\phi})\oplus (E'_{\bar \phi},\theta'_{\bar
\phi}),
$$
where $\calO_{LK}$ acts on the direct summand $E'_{\phi}$ via the
character $\phi$ and similarly for its bar counterpart.

\begin{proposition}\label{SOV correspondence}
Let $C_{\calM_0}: MIC(\calM_0)\to HIG(\calM'_0)$ be the Cartier
transform. Then for each $\phi\in \Phi$, one has $C_{\calM_0}(\calH^
1_{dR,\phi},\nabla_{\phi})=(E'_{\phi},\theta'_\phi)$. The same is
true for its bar counterpart.
\end{proposition}
\begin{proof}
Take $\lambda\in \calO_{LK}$. Applying Theorem 3.8 in \cite{OV} to
the morphism $\varsigma(\lambda): \calX_0\to \calX_0$ and the
objects $(\calO_{\calX_0},0)\in MIC(\calX_0),\
(\calO_{\calX_0'},0)\in HIG(\calX_0')$, one obtains the following commutative
diagram:
$$
  \xymatrix{
     (\calO_{\calX_0},0)\ar[r]^{C_{\calX_0}}\ar[d]^{\cong}_{\varsigma(\lambda)^{DR}_{*}}&  (\calO_{\calX_0'},0)  \ar[d]_{\cong}^{\varsigma'(\lambda)^{HIG}_{*}}   \\
     (\calO_{\calX_0},0)\ar[r]^{C_{\calX_0}} & (\calO_{\calX_0'},0).}
 $$
Applying the same theorem further to $f_0: \calX_0\to \calM_0$ one
obtains the second commutative diagram:
$$
  \xymatrix{
     (\calO_{\calX_0},0)\ar[r]^{C_{\calX_0}}\ar[d]^{ }_{R^1f_{0*}^{DR}}&  (\calO_{\calX_0'},0)  \ar[d]_{}^{R^1f_{0*}^{'HIG}}   \\
     (\calH^1_{dR},\nabla)\ar[r]^{C_{\calM_0}} & (E',\theta').}
 $$
It is clear that the isomorphism $\varsigma(\lambda)^{DR}_{*}$
(resp. $\varsigma'(\lambda)^{HIG}_{*}$) induces via the direct
image functor $R^1f_{0*}$ the isomorphism $\varsigma(\lambda)^{DR}_{*}:
(\calH^1_{dR},\nabla)\to (\calH^1_{dR},\nabla)$ (resp. the isomorphism
$\varsigma'(\lambda)^{HIG}_{*}:(E',\theta')\to (E',\theta')$). The
above two commutative diagrams yield the following commutative
diagram:
$$
  \xymatrix{
     (\calH^1_{dR},\nabla)\ar[r]^{C_{\calM_0}}\ar[d]^{\cong }_{\varsigma(\lambda)^{DR}_{*}}&  (E',\theta')  \ar[d]_{\cong}^{\varsigma'(\lambda)^{HIG}_{*}}   \\
     (\calH^1_{dR},\nabla)\ar[r]^{C_{\calM_0}} & (E',\theta').}
$$
Since the action of $\lambda\in \calO_{LK}$ on $(\calH^1_{dR},\nabla)$
(resp. $(E',\theta')$) is given by $\varsigma(\lambda)^{DR}_*$
(resp. $\varsigma'(\lambda)^{HIG}_{*}$), one has then
\begin{eqnarray*}
% \nonumber to remove numbering (before each equation)
   \varsigma'(\lambda)(C_{\calM_0}(\calH^1_{dR,\phi},\nabla_{\phi})) &=& C_{\calM_0}(\varsigma(\lambda)(\calH^1_{dR,\phi},\nabla_{\phi}))  \\
    &=&   C_{\calM_0}(\phi(\lambda)(\calH^1_{dR,\phi},\nabla_{\phi}))\\
    &=& \phi(\lambda)  C_{\calM_0}(\calH^1_{dR,\phi},\nabla_{\phi}),
\end{eqnarray*}
which implies that
$C_{\calM_0}(\calH^1_{dR,\phi},\nabla_{\phi})=(E'_{\phi},\theta'_{\phi})$.

\end{proof}

Let $0=F^2_{con}\subset F^1_{con}\subset F^0_{con}=\calH^1_{dR}$ be
the conjugate filtration of $\calH^1_{dR}$, which is flat with
respect to the Gauss-Manin connection (see \S 3 in \cite{Ka0}). For
a subbundle $W\subset \calH^1_{dR}$ we put
$Gr_{F_{con}}(W)=\bigoplus_{q=0}^{1} \frac{W\cap F_{con}^{q}}{W\cap
F_{con}^{q+1}}$. The $p$-curvature $\psi_{\nabla}$ of $\nabla$
defines the $F$-Higgs bundle
$$
\psi_{\nabla}:Gr_{F_{con}}(\calH^1_{dR})\to
Gr_{F_{con}}(\calH^1_{dR})\otimes F_{\calM_0}^*\Omega_{\calM_0}.
$$
As a reminder to the reader, we recall the definition of the
$F$-Higgs bundle: an $F$-Higgs bundle over a base $C$, which is
defined over $\bbF$, is a pair $(E,\theta)$ where $E$ is a vector
bundle over $C$, and $\theta$ is a bundle morphism $E\to E\otimes
F_C^*\Omega_C$ with the integral property $\theta\wedge\theta=0$.
The following lemma is a simple consequence of Katz's $p$-curvature
formula (Theorem 3.2, \cite{Ka4}), and it is also true in a general
context.
\begin{lemma}\label{grading of conjugate filtration}
Let $W\subset \calH^1_{dR}$ be a subbundle preserved by the
Gauss-Manin connection $\nabla$. Then the $F$-Higgs subbundle
$(Gr_{F_{con}}(W),\psi_{\nabla}|_{Gr_{F_{con}}(W)})$ defines a Higgs
subbundle of $(E',\theta')$ by the Cartier descent.
\end{lemma}
\begin{proof}
Since $\nabla$ preserves $F_{con}$, it induces the connection
$Gr_{F_{con}}\nabla$ on $Gr_{F_{con}}(\calH^1_{dR})$ by taking
grading. The operation of the $p$-curvature on this connection
commutes with taking the grading. It follows that
$\psi_{Gr_{F_{con}}\nabla}=Gr_{F_{con}}\psi_\nabla$, which is the
zero map. The relative Cartier isomorphism gives the isomorphism
$(Gr_{F_{con}}(\calH^1_{dR}),Gr_{F_{con}}\nabla)\cong
(\calF_{\calM_0}^*E',\nabla^{can})$, where $\nabla^{can}$ is the
canonical connection by the Cartier descent. Taking this isomorphism
for granted, we see that the inclusion $Gr_{F_{con}}(W)\subset
Gr_{F_{con}}(\calH^1_{dR})$ descends to the inclusion $F\subset E'$.
In other words, $\calF_{\calM_0}^*F$ is isomorphic to
$Gr_{F_{con}}(W)$ via the relative inverse Cartier isomorphism. By
Katz's formula (see \cite{Ka4}), the $F$-Higgs bundle
$(Gr_{F_{con}}(\calH^1_{dR}),\psi_{\nabla})$ descends to the Higgs
bundle $(E',\theta')$. Thus the $F$-Higgs subbundle
$(Gr_{F_{con}}(W),\psi_{\nabla}|_{Gr_{F_{con}}(W)})$ descends to the
Higgs subbundle $(F,\theta'|_{F})$ of $(E',\theta')$.
\end{proof}
For simplicity the above resulting Higgs subbundle is called the
\emph{Cartier descent} of $Gr_{F_{con}}(W,\nabla|_W)$.
\begin{theorem}\label{Cartier descent}
For each $\phi\in \Phi$, the following statements hold:
\begin{enumerate}
  \item
  $C_{\calM_0}(\calH^1_{dR,\sigma\phi},\nabla_{\sigma\phi})=\pi^*(E_{\phi},\theta_{\phi})$.
  \item
$C_{\calM_0}(\calH^1_{dR,\phi},\nabla_{\phi})$ is the Cartier
descent of $Gr_{F_{con}}(\calH^1_{dR,\phi},\nabla_{\phi})$.
\end{enumerate}
The similar statements hold for the bar counterpart.
\end{theorem}
\begin{proof}
By Proposition \ref{SOV correspondence} it is equivalent to show the
identification $(E_{\phi})'=E_{\sigma \phi}$. For this we take $\lambda\in
\calO_{LK}$. Then
$$\varsigma'(\lambda)((E_{\phi})')=(\varsigma(\lambda)(E_{\phi}))'=(\phi(\lambda)E_{\phi})'=\phi(\lambda)^p(E_{\phi})'=(\sigma\phi)(\lambda)(E_{\phi})'.$$
This proves the identification and therefore the first part of the theorem.
Because $\theta'|_{E'_{\phi}}=\theta'_{\phi}$, the second part is a
consequence of Proposition \ref{SOV correspondence} and Lemma
\ref{grading of conjugate filtration}.
\end{proof}

By this theorem we have a better understanding of the Higgs subbundles
in Proposition \ref{sss} (ii) under the Frobenius pull-backs. Now
assume that $d=[F_{\mathfrak p}:\bbQ_p]>1$ and write
$Hom_{\bbQp}(L_{\frakp},\bbQpp)=\{\phi_1,\cdots,\phi_d,\phi_1^*,\cdots,\phi_d^*\}$
as in Proposition \ref{isomorphisms between direct summands of
de Rham bundle}. We deduce the following
\begin{corollary}\label{HN filtration vs Hodge filtration}
For $2\leq i\leq d$, one has $
F_{\calM_0}^{*d-i+1}E_{\phi_i}=\calH^1_{dR,\phi_1^*}$. Moreover
under this identification, the Harder-Narasimhan filtration on
$F_{\calM_0}^{*d-i+1}E_{\phi_i}$ coincides with the Hodge filtration
on $\calH^1_{dR,\phi_1^*}$, which is induced from $\calH^1_{dR}$.
\end{corollary}
\begin{proof}
We look at the $i=d$ case first. Since $((E_{\phi_d})',0)\in HIG(\calM_0')$
is of trivial Higgs field,
$$C^{-1}_{\calM_0}((E_{\phi_d})',0)=(\calF_{\calM_0}^*(E_{\phi_d})',\nabla^{can}).$$
Because $F_{\calM_0}=\pi\circ \calF_{\calM_0}$, we have
$\calF_{\calM_0}^*(E_{\phi_d})'=F_{\calM_0}^*E_{\phi_d}$. On the
other hand, the above theorem says that
$$
C^{-1}_{\calM_0}((E_{\phi_d})',0)=C^{-1}_{\calM_0}(E'_{\sigma
\phi_d},0)= C^{-1}_{\calM_0}(E'_{\phi_
1^*},0)=(\calH^1_{dR,\phi_1^*},\nabla_{\phi_1^*}).
$$
It follows that
$(F_{\calM_0}^*E_{\phi_d},\nabla^{can})=(\calH^1_{dR,\phi_1^*},\nabla_{\phi_1^*})$.
Forgetting the connections, we obtain the result for $i=d$. For
$2\leq i\leq d-1$, the same argument for $E_{\phi_i}$ implies that
$F_{\calM_0}^*E_{\phi_{i}}=\calH^1_{dR,\phi_{i+1}}$. As $i+1\leq d$,
one must have $\calH^1_{dR,\phi_{i+1}}=E_{\phi_{i+1}}$. By iterating
the above arguments, one obtains the results for all $i$'s from
the $i=d$ case. From the second part of the proof of Proposition
\ref{sss}, we have known that $F_{\calM_0}^{*d-i+1}E_{\phi_i}$ is
not semistable. It suffices to show that the Hodge filtration on
$\calH^1_{dR,\phi_1^*}=F_{\calM_0}^*E_{\phi_d}$ is the
Harder-Narasimhan filtration. Since the Harder-Narasimhan filtration
is unique, it is equivalent to show that the sub line bundle
$E^{1,0}_{\phi_1^*}\subset \calH^1_{dR,\phi_1^*}$ is of maximal
degree. By Theorem \ref{instability}, the maximal degree does not
exceed $g-1$, which is exactly the degree of $E^{1,0}_{\phi_1^*}$.
This completes the proof.
\end{proof}
Taking the grading of
$(F_{\calM_0}^{*d-i+1}E_{\phi_i},\nabla^{can})$ with respect to the
Harder-Narasimhan filtration, one obtains a Higgs bundle different
from the original Higgs subbundle $(E_{\phi_i},0)$ of unitary type.
By the above result, the new Higgs bundle is exactly one of the
Higgs subbundles of unformizing type. Such a phenomenon might be
quite general for Higgs bundles over compact Shimura varieties in
char $p$.

\end{document}